\documentclass[11pt,reqno]{amsart}
\usepackage{amssymb,amsthm,amsmath,verbatim,enumerate,ifthen}
\usepackage[mathscr]{eucal}
\usepackage[utf8]{inputenc}
\usepackage[T1]{fontenc}
\def\N{\mathbb{N}}
\def\R{\mathbb{R}}
\def\CR{\mathscr{R}}
\def\S{\mathscr{S}}
\def\Q{\mathscr{Q}}

\def\A{\mathscr{A}}
\def\C{\mathscr{C}}

\def\B{\mathscr{E}}
\def\E{\mathscr{E}}

\def\T{\mathscr{T}}
\def\W{\mathscr{W}}
\def\TW{\mathscr{TW}}

\def\cl{\mathop{\mbox{\rm cl}}}

\def\qconv{\mathop{\mbox{\rm qconv}}\nolimits}
\def\dom{\mathop{\mbox{\rm dom}}\nolimits}
\def\codom{\mathop{\mbox{\rm codom}}\nolimits}
\def\epi{\mathop{\mbox{\rm epi}}\nolimits}
\def\graph{\mathop{\mbox{\rm graph}}\nolimits}
\def\div{\mathop{\text{\rm div}}\nolimits}

\def\Cdot{\!\cdot\!}

\newtheorem{theorem}{Theorem}[section]
\newtheorem*{theorem*}{Theorem}
\def\Thm#1#2{\ifthenelse{\equal{#1}{*}}{\begin{theorem*}#2\end{theorem*}}
             {\begin{theorem}\label{T#1}#2\end{theorem}}}
\newtheorem{Atheorem}{Theorem}

\def\THM#1#2{\begin{Atheorem}\label{T#1}#2\end{Atheorem}}
\def\thm#1{Theorem~\ref{T#1}}

\newtheorem{proposition}[theorem]{Proposition}
\newtheorem*{proposition*}{Proposition}
\def\Prp#1#2{\ifthenelse{\equal{#1}{*}}{\begin{proposition*}#2\end{proposition*}}
             {\begin{proposition}\label{P#1}#2\end{proposition}}}
\def\prp#1{Proposition~\ref{P#1}}

\newtheorem{corollary}[theorem]{Corollary}
\newtheorem*{corollary*}{Corollary}
\def\Cor#1#2{\ifthenelse{\equal{#1}{*}}{\begin{corollary*}#2\end{corollary*}}
             {\begin{corollary}\label{C#1}#2\end{corollary}}}
\def\cor#1{Corollary~\ref{C#1}}

\newtheorem{lemma}[theorem]{Lemma}
\newtheorem*{lemma*}{Lemma}
\def\Lem#1#2{\ifthenelse{\equal{#1}{*}}{\begin{lemma*}#2\end{lemma*}}
             {\begin{lemma}\label{L#1}#2\end{lemma}}}
\def\lem#1{Lemma~\ref{L#1}}

\theoremstyle{definition}
\newtheorem{remark}[theorem]{Remark}
\newtheorem*{remark*}{Remark}
\def\Rem#1#2{\ifthenelse{\equal{#1}{*}}{\begin{remark*}\rm #2\end{remark*}}
             {\begin{remark}\label{R#1}\rm #2\end{remark}}}

\newtheorem{example}[theorem]{Example}
\newtheorem*{example*}{Example}
\def\Exa#1#2{\ifthenelse{\equal{#1}{*}}{\begin{example*}\rm #2\end{example*}}
             {\begin{example}\label{Ex#1}\rm #2\end{example}}}

\def\eq#1{{\rm(\ref{E#1})}}
\def\Eq#1#2{\ifthenelse{\equal{#1}{*}}
  {\begin{equation*}\begin{aligned}[]#2\end{aligned}\end{equation*}}
  {\begin{equation}\begin{aligned}[]\label{E#1}#2\end{aligned}\end{equation}}}

\begin{document}
\vspace{5mm}

\date{\today}

\title[]{Convexity properties of functions defined on metric Abelian groups}

\author[W. Fechner]{W{\l}odzimierz Fechner}
\address{Institute of Mathematics, Lodz University of Technology, ul. W\'olcza\'nska 215, 90-924 \L\'od\'z, Poland}
\email{wlodzimierz.fechner@p.lodz.pl}

\author[Zs. P\'ales]{Zsolt P\'ales}
\address{Institute of Mathematics, University of Debrecen, 
H-4002 Debrecen, Pf.\ 400, Hungary}
\email{pales@science.unideb.hu}


\thanks{The research of the second author was supported by the EFOP-3.6.1-16-2016-00022 and the EFOP-3.6.2-16-2017-00015 projects. These projects are co-financed by the European Union and the European Social Fund.}

\begin{abstract}
The notions of quasiconvexity, Wright convexity and convexity for functions defined on a metric Abelian group are introduced. Various characterizations of such functions, the structural properties of the functions classes so obtained are established and several well-known results are extended to this new setting.
\end{abstract}

\maketitle

\section{Introduction}

Let $X$ be a linear space and $t\in[0,1]$. 
A subset $D\subseteq X$ is termed \emph{$t$-convex} if, for all $x,y\in D$, 
\Eq{*}{
  tx+(1-t)y\in D.
}
Analogously, a function $f:D\to\R$ is called \emph{$t$-quasiconvex, $t$-Wright convex, and $t$-convex} if $D$ is a $t$-convex set and, for all $x,y\in D$, the respective inequality
\Eq{tC}{
  f(tx+(1-t)y)&\leq \max(f(x),f(y)),\\
  f(tx+(1-t)y)+f((1-t)x+ty)&\leq f(x)+f(y),\\
  f(tx+(1-t)y)&\leq tf(x)+(1-t)f(y),
}
holds. The $\frac12$-convex sets are said to be \emph{midpoint convex} and the $\frac12$-convex functions are usually called \emph{Jensen convex}. The structure and properties of $t$-convex sets and $t$-quasiconvex, $t$-Wright convex, and $t$-convex functions and their generalizations have been investigated in a large number of recent papers, see e.g. \cite{Ada03b,GilNikPal04,GilPal08,Kom92,Kom03,Kuh84,Lew07,Lew08,Lew09,MakNikPal91,MakPal09b,MakPal10a,MakPal11,Mat92c,MatPyc95b,MatWro96,Ng87b,NikPal01,Olb11b,Olb11a,Olb13,Olb15d,Olb15b,Olb17b,OlbPal18,Pyc06,Pal13,Wri54}.

As a consequence of a result by Dar\'oczy and P\'ales \cite{DarPal87}, every $t$-convex function (where $t\in\,]0,1[\,$) is automatically Jensen convex and hence $\mathbb{Q}$-convex, i.e., it is $r$-convex for all rational numbers $r\in[0,1]$ (cf.\ \cite{Kuc85}). The following more general result about $t$-convexity was established by Kuhn \cite{Kuh84}.

\THM{Kuhn}{If $D$ contains at least two points and $f:D\to\R$ is a $t$-convex function for some $t\in\,]0,1[\,$, then $f$ is $s$-convex for all $s\in\Q(t)\cap[0,1]$, where $\Q(t)$ denotes the smallest subfield of $\R$ containing $t$. Furthermore, for every subfield $F$ of $\R$, there exists a function $f:D\to\R$ which is $t$-convex if and only if $t\in F\cap[0,1]$.}

The following result is the multivariable extension of the $t$-convexity property. 

\THM{Hconv}{Let $F$ be a subfield of $\R$ and $f:D\to\R$ be $t$-convex for all $t\in F\cap[0,1]$. Then, for all $n\in\N$, $x_1,\dots,x_n\in D$, $t_1,\dots,t_n\in F\cap[0,1]$ with $t_1+\dots+t_n=1$, the following inequality holds:
\Eq{*}{
  f\bigg(\sum_{i=1}^nt_ix_i\bigg)
  \leq \sum_{i=1}^n t_if(x_i).
}}

Another classical theorem is due to Bernstein and Doetsch \cite{BerDoe15} (see also \cite{Kuc85}).

\THM{BDT}{Let $D$ be an open convex subset of a normed linear space and let $f:D\to\R$ be a Jensen convex function which is bounded from above on a nonvoid open subset of $D$. Then $f$ is continuous and convex, that is, $t$-convex for all $t\in[0,1]$.}

In the paper \cite{MakNikPal91} the question whether $t$-Wright convexity implies Jensen convexity was investigated and an affirmative answer was proved if $t$ is a rational number. It was also shown that, for a transcendental $t$, this implication is not true. Furthermore, it turned out that for some second degree algebraic numbers the answer is positive whereas for some second degree algebraic numbers is negative. Bernstein--Doetsch-type theorems for Wright convex functions were established by Olbry\'s \cite{Olb13} and by Lewicki \cite{Lew08,Lew09}. On the other hand, in \cite{GilNikPal04} Bernstein--Doetsch-type theorems were proven for quasiconvex functions.

All the above mentioned results motivate to investigate the analogous problems in a more general setting. In our previous paper \cite{FecPal19a} we have defined the convexity of sets in metric Abelian groups with the help of endomorphisms. The purpose of this paper is to adopt and extend this definition to functions and therefore to investigate the associated notions of quasiconvexity, Wright convexity and convexity. Some of our results will generalize \thm{Kuhn} and \thm{Hconv} and also the above described statements.

\section{Metric Abelian groups and convexity of subsets}

In this section we briefly recall the terminology, the notations and all the results from \cite{FecPal19a} which will be instrumental for our approach.

Let $(X,+)$ be an Abelian group and let $\E(X)$ denote the family of all endomorphisms. Then $(\E(X),+,\circ)$ is a ring. Thus, every $T\in\E(X)$ generates an endomorphism $\widetilde{T}:\E(X)\to\E(X)$ defined by $\widetilde{T}(S):=T\circ S$. For a family $\T\subseteq \E(X)$ we denote $\widetilde{\T}:=\{\widetilde{T} \mid T \in \T\}$. Finally, $I$ stands for the identity map of $X$. The multiplication of the elements of $X$ by natural numbers is introduced by 
\Eq{*}{
  1\Cdot x:=x,\qquad\mbox{and}\qquad (n+1)\Cdot x:=n\Cdot x+x \qquad(x\in X,\,n\in\N).
}
The mapping $\pi_n(x):=n\Cdot x$ is always an endomorphism of $X$. We say that $(X,+)$ is \emph{divisible by} $n\in\N$ if the map $\pi_n$ is a bijection (and hence an automorphism) of $X$. In this case, for $x\in X$, the element $\pi_n^{-1}(x)$ is denoted as $\frac{1}{n}\Cdot x$. The set of natural numbers $n$ for which $X$ is uniquely divisible by $n$ is a multiplicative subsemigroup of $\N$ whose unit element is $1$, will be denoted by $\div(X)$.

For a subset $A\subseteq X$ and $n\in\N$, we say that $A$ is \emph{$n$-convex} if 
\Eq{*}{
  \{n\Cdot x\mid x\in A\}
  =\{x_1+\dots+x_n\mid x_1,\dots,x_n\in A\}.
}
For properties of $n$-convex sets, we refer to the paper \cite{JarPal15}. In particular, by \cite[Proposition 2]{JarPal15}, we have that if a set is $n$- and $m$-convex, then it is also $(nm)$-convex.

In the case when $(X,+)$ is equipped with a \emph{translation invariant metric $d$}, we say that $(X,+,d)$ is a \emph{metric Abelian group}. Metric groups are automatically topological groups in which the \emph{$d$-norm} $\|\Cdot\|_d:X\to\R$ is defined as $\|x\|_d:=d(x,0)$. The subadditivity of $\|\Cdot\|_d$ implies that $\|n\Cdot x\|_d\leq n\|x\|_d$ for all $x\in X$ and $n\in\N$. The equality here, may not be valid.  

An endomorphism $T:X\to X$ is called \emph{$d$-bounded} if there exists $c\geq0$ such that $\|T(x)\|_d\leq c\|x\|_d$ for all $x\in X$. The smallest number $c$ satisfying this condition is called the \emph{$d$-norm} of $T$ and is denoted by $\|T\|^*_d$. The symbol $\E^d(X)$ will denote the subring of $\E(X)$ of all $d$-bounded endomorphisms. More generally, for $\T\subseteq\E(X)$, the symbol $\T^d$ denotes the $d$-bounded elements of $\T$. The smallest number $c$ such that 
$\|n\Cdot x\|_d\leq c\|x\|_d$ for all $x\in X$, that is $\|\pi_n\|^*_d$, will simply be denoted by $\|n\|^*_d$.

For $n\in\N$, the measure of injectivity of the map $\pi_n$ is the largest number $\mu_d(n)$ such that
\Eq{mun}{
   \mu_d(n)\|x\|_d\leq \|n\Cdot x\|_d \qquad(x\in X).
}
Using these notations, we can now formulate an extension of the celebrated R{\aa}dstr\"{o}m Cancellation Theorem (cf.\cite{Rad52b}) which we proved in \cite{FecPal19a}.

\Thm{RCT}{Let $(X,+,d)$ be a metric Abelian group and let $n_0\in\N$ such that $\mu_d(n_0)>1$. Let $A\subseteq X$ be an arbitrary subset, let $B\subseteq X$ be closed and $n_0$-convex subset, and $C\subseteq X$ be a $d$-bounded nonempty subset such that $A+C\subseteq B+C$. Then $A\subseteq B$.}

The \emph{$d$-spectral radius of an endomorphism $T\in\E^d(X)$} is defined as
\Eq{*}{
  \rho_d(T):=\limsup_{m\to\infty}\sqrt[m]{\|T^m\|^*_d}.
}
The following result is a generalization of the so-called Neumann invertibility theorem.

\Thm{NIT}{Let $(X,+,d)$ be a complete metric Abelian group and let $T\in\E^d(X)$ such that $\rho_d(T)<1$. Then $I-T$ is an invertible element of $\E^d(X)$, furthermore,
\Eq{*}{
  (I-T)^{-1}=\sum_{k=0}^\infty T^k.
}}

Given an endomorphism $T\in\E(X)$, we say that a subset $D\subseteq X$ is \emph{$T$-convex} if, for all $x,y\in D$, 
\Eq{*}{
  T(x)+(I-T)(y)\in D.
}
This condition is equivalent to the inclusion
\Eq{*}{
  T(D)+(I-T)(D)\subseteq D.
}
If $\T\subseteq\E(X)$, then a set $D\subseteq X$ is called \emph{$\T$-convex} if it is $T$-convex for all $T\in\T$. The class of $\T$-convex subsets of $X$ is denoted by $\C_\T(X)$ in what follows. In the particular case when $(X,+)$ is the additive group of a vector space and $T=tI$ for some $t\in[0,1]$, instead of $T$-convexity, we briefly speak about \emph{$t$-convexity} which is a commonly accepted notion (cf.\ \cite{Kuh84}). If $X$ is a uniquely $2$-divisible Abelian group, and $T=\frac12\Cdot I$, that is, $T(x):=\frac12\Cdot x$, then $T$-convex sets will also be termed \emph{midpoint convex}. One can immediately see that if the group $X$ is divisible by some $n\in \N$ and $T=\frac{1}{n}\Cdot I$, then $T$-convexity is equivalent to $n$-convexity defined in the previous section. It is obvious but useful to observe that a subset $D\subseteq X$ is $T$-convex if and only if, for all $p\in D$, 
\Eq{Tp}{
  T(D-p)\subseteq D-p.
}

Now, given a nonempty subset $D\subseteq X$, we consider the collection of endomorphisms $T$ of $X$ that make $D$ to be $T$-convex: 
\Eq{*}{
  \T_D:=\{T\in\E(X)\mid \mbox{$D$ is $T$-convex}\}.
}
It is obvious that, for every set $D$, we have $0,I\in\T_D$ and $0,I\in\T^d_D$ (if $X$ is a metric Abelian group). The next result describes a convexity property of $\T_D$.

\Thm{P1}{Let $D\subseteq X$ be a nonempty set. Then $\T_D$ is a $\widetilde\T_D$-convex subset of $\E(X)$. If $(X,+,d)$ is a metric Abelian group, then $\T^d_D$ is a $\widetilde\T^d_D$-convex subset of $\E^d(X)$. In particular, these sets are closed with respect to the composition of maps.}

\Cor{1}{Let $D\subseteq X$ be a nonempty set. Then $\T_D$ and also $\T^d_D$ (if $(X,+,d)$ is a metric Abelian group) are closed under multiplication and under the mappings 
\Eq{maps}{
T\mapsto I-T\qquad\mbox{and}\qquad (T,S)\mapsto T\circ S+(I-T)\circ(I-S).
}}

In the next result, we provide conditions ensuring that $T$-convexity implies midpoint convexity. 

\Thm{2}{Let $(X,+,d)$ be a complete metric uniquely 2-divisible Abelian group and $T\in\E^d(X)$ such that $\rho_d(2\Cdot T-I)<1$. Then, for every nonempty $T$-convex set $D\subseteq X$, the set $\cl(\T^d_D)$ is a midpoint convex subset of $\E^d(X)$. Furthermore, every closed $T$-convex subset of $X$ is also midpoint convex.}

The following result will be instrumental when investigating $T$-Wright convex functions.

\Cor{nkc2}{Let $(X,+,d)$ be a metric  Abelian group, let $n_0\in\N$ such that $\mu_d(n_0)>1$ and let $D$ be a closed bounded $n_0$-convex set. Let $n\in\N$ and $T_1,\dots,T_n\in\T^d_D$ be such that $T:=T_1+\cdots+T_n$ is a bijection with $T^{-1}\in\E^d(X)$. Then, for all $k\in\{1,\dots,n-1\}$, we have $T^{-1}\circ(T_1+\cdots+T_k)\in\T_D^d$.}

\section{$T$-quasiconvex functions}

Assume that $(X,+)$ is an Abelian group.
Given an endomorphism $T\in\E(X)$, we say that a function $f:D\to[-\infty,+\infty[\,$ is \emph{$T$-quasiconvex} if $D$ is a $T$-convex subset of $X$ and, for all $x,y\in D$, 
\Eq{*}{
  f(T(x)+(I-T)(y))\leq\max(f(x),f(y)).
}
If $\T\subseteq\E(X)$, then a function $f:D\to[-\infty,+\infty[\,$ is called \emph{$\T$-quasiconvex} if, for all $T\in\T$, it is $T$-quasiconvex. In what follows, the class of $\T$-quasiconvex functions defined on $D$ is denoted by $\Q_\T(D)$. If $(X,+)$ is the additive group of a linear space and $f$ is $tI$-quasiconvex for some $t\in[0,1]$, then we say that $f$ is \emph{$t$-quasiconvex}. If $f$ is $t$-quasiconvex for all $t\in[0,1]$, then it is called a \emph{quasiconvex function} (in the standard sense). If $(X,+)$ is a uniquely 2-divisible group, and $f$ is $\frac12\Cdot I$-quasiconvex, then it is called \emph{midpoint quasiconvex} (cf.\ \cite{NikPal98}).  Recall that the characteristic function of a set $S\subseteq X$ is defined by
\Eq{*}{
  \chi_S(x):=\begin{cases}
           1 &\mbox{if }x\in S,\\
           0 &\mbox{if }x\in X\setminus S.
           \end{cases}
}

\Prp{LS}{Let $\T\subseteq\E(X)$ be a nonempty subset and let $D$ be a $\T$-convex set. Then a function $f:D\to[-\infty,+\infty[\,$ is $\T$-quasiconvex if and only if, for all $c\in[-\infty,+\infty[\,$, the level sets
\Eq{*}{
   D^c_f:=\{x\in D\colon f(x)\leq c\}
}
are $\T$-convex. On the other hand, a set $S$ is $\T$-convex if and only if the negative of its characteristic function $\chi_S$ is $\T$-quasiconvex.}

\begin{proof}
Assume that $T\in \T$, the function $f:D\to[-\infty,+\infty[\,$ is $T$-quasiconvex and fix a $c\in[-\infty,+\infty[\,$ arbitrarily. We will show that the set $D^c_f$ is $T$-convex. We have to show that, for all $p \in D^c_f$, the inclusion $T(D^c_f - p) \subseteq D^c_f - p$ is true. We have $f(p) \leq c$ and let $y \in  T(D^c_f - p)$. There exists $x \in D^c_f$ such that $y=T(x-p)$. Clearly $f(x) \leq c$. By $T$-quasiconvexity of $f$ we have 
 \Eq{*}{
 f(T(x-p)+p)  =f(T(x) + (I-T)(p)) \leq  \max(f(x),f(p)) \leq c,
}
which gives us that $y+p = T(x-p)+p \in D^c_f$, which was to be shown.

To prove the converse implication fix $T\in \T$, $x, y \in D$ and assume that the level sets of a function $f:D\to[-\infty,+\infty[\,$ are $T$-convex sets. Without loss of generality we can assume that $f(y)\leq f(x)$. Then $y \in D^{f(y)}_f\subseteq D^{f(x)}_f$ and in particular the set $D^{f(x)}_f$ is $T$-convex. Therefore, we have $T(D^{f(x)}_f - y) \subseteq D^{f(x)}_f - y$. Consequently, $T(x-y) \in D^{f(x)}_f - y$, which means that 
\Eq{*}{f(T(x) + (I-T)(y)) = f(T(x-y) + y) \leq f(x) = \max(f(x),f(y)).
}
This proves the $T$-quasiconvexity of $f$ when level sets are $T$-convex.

Finally, assume that we are given a set $S\subseteq X$ and $T\in \T$. Note that $T$-quasiconvexity of $-\chi_S$ can be directly rewritten as follows:
\Eq{*}{
\chi_S(T(x) + (I-T)(y)) \geq \min (\chi_S(x),\chi_S(y)) 
}
for all $x, y \in X$. This inequality in turn  is equivalent to (use the definition of characteristic function):
\Eq{*}{
x, y \in S \, \Longrightarrow \, T(x) + (I-T)(y) \in S,
}
which is precisely the $T$-convexity of the set $S$.
\end{proof}

\Thm{TQ}{Let $\T\subseteq\E(X)$ be a nonempty subset. Then we have the following statements.
\begin{enumerate}[(i)]
 \item If $D$ is a $\T$-convex set, then $\Q_\T(D)$ contains negative of the characteristic functions of all $\T$-convex subsets of $D$ and, for every $c\in[-\infty,+\infty[\,$ and $f\in \Q_\T(D)$, we have $f+c\in\Q_\T(D)$.
 \item If $D$ is a $\T$-convex set, then $\Q_\T(D)$ is closed with respect to the pointwise supremum, the pointwise chain infimum, and the pointwise convergence.
 \item If $f:D\to[-\infty,+\infty[\,$ and $g:E\to[-\infty,+\infty[\,$ are $\T$-quasiconvex functions, then the function $f\diamond g:D+E\to[-\infty,+\infty[\,$ defined by
\Eq{*}{
  (f\diamond g)(x):=\inf\{\max(f(u),g(v)):u\in D,\,v\in E,\,u+v=x\}
}
is $\T$-quasiconvex on $D+E$.
\item If an endomorphism $A\in\E(X)$ commutes with any member of $\T$ and $f\in\Q_\T(D)$, then $f\circ A\in\Q_\T(A^{-1}(D))$ and $f\circ A^{-1}\in\Q_\T(A(D))$, where $f\circ A^{-1}:A(D)\to[-\infty,+\infty[\,$ is defined by
\Eq{fAi}{
  (f\circ A^{-1})(x):=\inf_{u\in A^{-1}(x)\cap D} f(u).
}
\end{enumerate}}

\begin{proof}

The proof of (i) is obvious.

Proof of (ii). To prove the $\T$-quasiconvexity of the pointwise supremum of a family $\{f_\alpha \mid \alpha \in I\}$ of  $\T$-quasiconvex functions defined on $D$ fix $T \in \T$, $x, y \in D$ and $\varepsilon >0$. Let  $f:D\to[-\infty,+\infty[\,$ be given by $f=\sup \{f_\alpha \mid \alpha \in I\}$. To show that $f$ is $T$-quasiconvex observe that there exists some $\alpha_0\in I$ such that 
\Eq{*}{
 f(T(x)+(I-T)(y)) \leq f_{\alpha_0}(T(x)+(I-T)(y)) + \varepsilon \leq \max(f_{\alpha_0}(x),f_{\alpha_0}(y)) + \varepsilon\leq \max(f(x),f(y)) + \varepsilon.
}
Since $\varepsilon >0$ is arbitrary small, then $f$ is $T$-quasiconvex.

To justify the $\T$-quasiconvexity of the pointwise infimum of a chain $\{f_\alpha \mid \alpha \in I\}$ of  $\T$-quasiconvex functions defined on $D$ fix $T \in \T$, $x, y \in D$ and $\varepsilon >0$. Let  $f:D\to[-\infty,+\infty[\,$ be given by $f=\inf \{f_\alpha \mid \alpha \in I\}$. To show that $f$ is $T$-quasiconvex observe that there exist some $\alpha_x, \alpha_y\in I$ such that $f_{\alpha_x}(x) \leq f(x) + \varepsilon$ and $f_{\alpha_y}(y) \leq f(y) + \varepsilon$. Since family $\{f_\alpha \mid \alpha \in I\}$ forms a  chain, there exists $\alpha_0\in \{\alpha_x,\alpha_y\}$ such that $f_{\alpha_0}=\min(f_{\alpha_x}, f_{\alpha_y})$ and then we have
\Eq{*}{
 f(T(x)+(I-T)(y)) &\leq f_{\alpha_0}(T(x)+(I-T)(y))  \leq \max(f_{\alpha_0}(x),f_{\alpha_0}(y)) \\&\leq \max(f_{\alpha_x}(x),f_{\alpha_y}(y)) \leq \max(f(x) + \varepsilon,f(y)+ \varepsilon) = \max(f(x),f(y)) + \varepsilon.
}
Again, since $\varepsilon >0$ is arbitrary small, we obtain that $f$ is $T$-quasiconvex.

To show the $\T$-quasiconvexity of the pointwise limit $(f_n)$ of  $\T$-quasiconvex functions defined on $D$ fix $T \in \T$ and $x, y \in D$. We have 
\Eq{*}{
 f(T(x)+(I-T)(y)) = \lim_{n \to +\infty}  f_n(T(x)+(I-T)(y)) \leq \lim_{n \to +\infty} \max(f_n(x),f_n(y)) = \max(f(x),f(y)).
}
The resulted equality proves that $f$ is $T$-quasiconvex.
 
Proof of (iii). For the $\T$-quasiconvexity of the function $f\diamond g$, let $x,y\in D+E$. We need to prove, for all $T\in\T$, that 
\Eq{f*g}{
  (f\diamond g)(T(x)+(I-T)(y))\leq\max((f\diamond g)(x),(f\diamond g)(y)).
}
Let $c\in \,](f\diamond g)(x),+\infty[\,$ and $d\in \,](f\diamond g)(y),+\infty[\,$ be arbitrary. Then, by the definition of $f\diamond g$, there exist $u,v\in D$ such that
\Eq{*}{
  \max(f(u),g(x-u))<c \qquad\mbox{and}\qquad \max(f(v),g(y-v))<d.
}
Let $T\in\T$ be fixed. Then, using the definition of $f\diamond g$ and the $\T$-quasiconvexity of $f$ and $g$, we obtain
\Eq{*}{
  (f\diamond g)(T(x)+(I-T)(y))
  &\leq \max\big(f(T(u)+(I-T)(v),g(T(x-u)+(I-T)(y-v))\big)\\
  &\leq \max\big(\max(f(u),f(v)),\max(g(x-u),g(y-v))\big) \\
  &= \max(f(u),g(x-u),f(v),g(y-v))<\max(c,d).
}
Upon taking the limits $c\searrow (f\diamond g)(x)$ and $d\searrow (f\diamond g)(y)$, the inequality \eq{f*g} follows.

Proof of (iv). To verify the $\T$-quasiconvexity of the function $f\circ A$, let $x,y\in A^{-1}(D)$ and let $T\in\T$ be fixed. Then $A(x),A(y) \in D$, hence the $T$-quasiconvexity of $f$ yields
\Eq{*}{
  (f\circ A)(T(x)+(I-T)(y))&=f(T(A(x))+(I-T)(A(y)) \\
  &\leq \max(f(A(x)),f(A(y)))=\max((f\circ A)(x),(f\circ A)(y)).
}
Finally, we show the $\T$-quasiconvexity of the function $f\circ A^{-1}$. For this proof, let $x,y\in A(D)$. For the proof of the inequality 
\Eq{fA+}{
  (f\circ A^{-1})(T(x)+(I-T)(y))\leq\max((f\circ A^{-1})(x),(f\circ A^{-1})(y))
}
choose $c\in \,](f\circ A^{-1})(x),+\infty[\,$ and $d\in \,](f\circ A^{-1})(y),+\infty[\,$ arbitrarily. Then, there exist $u,v\in D$ such that $A(u)=x$, $A(v)=y$ and $f(u)<c$, $f(v)<d$. Then, using $A(T(u)+(I-T)(v))=T(x)+(I-T)(y)$, and the $T$-quasiconvexity of $f$, we get
\Eq{*}{
  (f\circ A^{-1})(T(x)+(I-T)(y))
  \leq f(T(u)+(I-T)(v))\leq \max(f(u),f(v))<\max(c,d).
}
Now, upon taking the limits $c\searrow (f\circ A^{-1})(x)$ and $d\searrow (f\circ A^{-1})(y)$, the inequality \eq{fA+} follows.
\end{proof}

Using assertion (ii) of \thm{TQ}, it follows that, for every function $f:D\to[-\infty,+\infty[\,$ defined on a $\T$-convex set $D\subseteq X$, the function $\qconv_\T(f):D\to[-\infty,+\infty[\,$ defined as
\Eq{*}{
   \qconv_\T(f)(x):=\sup\{g(x)\mid g\in\Q_\T(D),\,g\leq f\} \quad(x\in D)
}
is the largest $\T$-quasiconvex function which is not greater than $f$ on $D$. This function will be called the \emph{$\T$-quasiconvex envelope of $f$}.

Now, given a function $f:D\to[-\infty,+\infty[\,$, we consider the collection of endomorphisms $T\in\E(X)$ that make $f$ to be $T$-quasiconvex: 
\Eq{*}{
  \T_f:=\{T\in\E(X)\mid \mbox{$f$ is $T$-quasiconvex}\}.
}
It is obvious that, for every function $f$, we have $0,I\in\T_f$ and $0,I\in\T^d_f$ (if $(X,+,d)$ is a metric Abelian group). The next result shows some structural properties of $\T_f$ and $\T_f^d$.

\Thm{P1f}{Let $D\subseteq X$ and let $f:D\to[-\infty,+\infty[\,$ be an arbitrary function. Then $\T_f$ is a $\widetilde\T_f$-convex subset of $\E(X)$. If $(X,+,d)$ is a metric Abelian group, then $\T^d_f$ is a $\widetilde\T^d_f$-convex subset of $\B(X)$.}

\begin{proof} Let $T,T_1,T_2\in\T_f$ and set $S:=T\circ T_1+(I-T)\circ T_2$. Then, $f$ is $T_1$-, $T_2$- and $T$-quasiconvex, therefore, for all $x,y\in D$, we have
\Eq{*}{
  f(T_1(x)+(I-T_1)(y))\leq\max(f(x),f(y)), \qquad f(T_2(x)+(I-T_2)(y))\leq\max(f(x),f(y)).
}
Consequently,
\Eq{*}{
  f(S(x)+(I-S)(y))
  &=f\big(T(T_1(x)+(I-T_1)(y))+(I-T)(T_2(x)+(I-T_2)(y))\big)\\
  &\leq \max(f\big(T_1(x)+(I-T_1)(y)\big),f\big(T_2(x)+(I-T_2)(y)\big))
  \leq \max(f(x),f(y)).
}
This means that $f$ is $S$-quasiconvex, hence $S\in\T_f$.
This yields that $\T_f$ is $\widetilde{T}$-convex for all $T\in\T_f$, which was to be proved.

The proof of the second assertion is completely analogous. 
\end{proof}

The next result follows from \thm{P1f} exactly in the same manner as \cor{1} was deduced from \thm{P1}.

\Cor{1f}{Let $D\subseteq X$ and $f:D\to[-\infty,+\infty[\,$. Then $\T_f$ and also $\T^d_f$ (if $(X,+,d)$ is a metric Abelian group) is closed under multiplication and under the mappings in \eq{maps}.}

In the following statement we show that $T$-quasiconvexity implies the midpoint quasiconvexity under certain conditions on $T$, $f$, and $X$.

\Thm{2f}{Let $(X,+,d)$ be a uniquely 2-divisible metric Abelian group and $T\in\E^d(X)$ such that $\rho_d(2\Cdot T-I)<1$ and let $D$ be a closed $T$-convex set. Then, for every function $f\in\Q_T(D)$, the set $\cl(\T^d_f)$ is a midpoint convex subset of $\E^d(X)$. Furthermore, every lower semicontinuous function $f\in\Q_T(D)$ is also midpoint quasiconvex on $D$.}

\begin{proof} In view of \thm{2}, we have that $D$ is a midpoint convex set. Let $f\in\Q_T(D)$ and define the sequence of endomorphisms $T_n$ by 
\Eq{*}{
  T_n=\frac{1}{2}\Cdot \Big(I+(2\Cdot T-I)^{2^{n-1}}\Big).
}
By induction, one can see that this sequence satisfies the recursion 
\Eq{*}{
  T_1:=T,\qquad T_{n+1}:=T_n^2+(I-T_n)^2 \quad(n\in\N).
}
Then, by the last assertion of \cor{1f}, it follows that $T_n\in\T^d_f$ for all $n\in\N$. 

The condition $\rho_d(2\Cdot T-I)<1$ implies that $T_n$ converges to $\frac12\Cdot I$. If $R,S\in\cl(\T^d_f)$, then there exist sequences $R_n$ and $S_n$ in $\T^d_f$ converging to $R$ and $S$, respectively. By \thm{P1f}, for all $n\in\N$, we have that $T_n\circ R_n+(I-T_n)\circ S_n\in\T^d_f$. Upon taking the limit, it follows that $\frac12\Cdot (R+S)\in\cl(\T^d_f)$. This implies that $\cl(\T^d_f)$ is  a midpoint convex set.

To complete the proof, assume that $f$ is also a lower semicontinuous function. To prove its midpoint quasiconvexity, let $x,y\in D$. Then the midpoint convexity of the set $\cl(\T^d_f)$ and $0,I\in\cl(\T^d_f)$ imply that $\frac12\Cdot I+\frac12\Cdot 0=\frac12\Cdot I\in\cl(\T^d_f)$. Therefore, there exists a sequence of operators $S_n\in\T^d_f$ which converges to $\frac12\Cdot I$. Thus, for all $n\in\N$, 
\Eq{*}{
  f(S_n(x)+(I-S_n)(y))\leq\max(f(x),f(y)).
}
Upon taking the limit $n\to\infty$ and using the lower semicontinuity of $f$, it follows that 
\Eq{*}{
  f\big(\tfrac12\Cdot(x+y)\big)\leq\max(f(x),f(y)).
}
Therefore, $f$ is midpoint quasiconvex on $D$.
\end{proof}

The following result presents a further invariance property of $\T^d_f$.

\Thm{qc}{Assume that $(X,+,d)$ is a metric Abelian group, $n_0\in\N$ is such that $\mu_d(n_0)>1$ and $D$ is a closed set. Let $f:D\to[-\infty,\infty[\,$ be a lower semicontinuous function whose level sets $D_f^c$ are bounded $n_0$-convex for all $c\in\R$. Let $n\in\N$ and $T_1,\dots,T_n\in\T^d_f$ be such that $T:=T_1+\cdots+T_n$ is a bijection with $T^{-1}\in\E^d(X)$. Then, for all $k\in\{1,\dots,n-1\}$, we have $T^{-1}\circ(T_1+\cdots+T_k)\in\T^d_f$.}

\begin{proof} The lower semicontinuity of $f$ implies that the level sets 
$D_f^c$ are closed bounded $n_0$-convex subsets of the closed set $D$ for all $c\in\R$. In view of \prp{LS}, it follows that these level sets are $T_1$-, \dots, $T_n$-convex. Now, applying \cor{nkc2}, we obtain that all these level sets are $\big(T^{-1}\circ(T_1+\cdots+T_k)\big)$-convex. Hence, again by \prp{LS}, we get that $f$ is $\big(T^{-1}\circ(T_1+\cdots+T_k)\big)$-quasiconvex.
\end{proof}

\section{$T$-Wright convex and $T$-Wright affine functions}

Assume that $(X,+)$ is an Abelian group.
For an endomorphism $T\in\E(X)$, we say that a function $f:D\to[-\infty,+\infty[\,$ is \emph{$T$-Wright convex} if $D$ is a $T$-convex subset of $X$ and, for all $x,y\in D$, 
\Eq{*}{
  f(T(x)+(I-T)(y))+f((I-T)(x)+T(y))\leq f(x)+f(y).
}
If $\T\subseteq\E(X)$, then a function $f:D\to[-\infty,+\infty[\,$ is called \emph{$\T$-Wright convex} if, for all $T\in\T$, it is $T$-Wright convex. The class of $\T$-Wright convex functions defined on $D$ is denoted by $\W_\T(D)$. If $(X,+)$ is the additive group of a linear space and $f$ is $t\Cdot I$-Wright convex for some $t\in[0,1]$, then we say that $f$ is \emph{$t$-Wright convex}. If $f$ is $t$-Wright convex for all $t\in[0,1]$, then it is called a \emph{Wright convex} function (in the standard sense) (cf.\ \cite{Wri54}). If $(X,+)$ is a uniquely 2-divisible group, then $\frac12\Cdot I$-Wright convex functions are called \emph{Jensen convex} (cf.\ \cite{Kuc85}). 

\Thm{TW}{Let $\T\subseteq\E(X)$ be a nonempty subset. Then we have the following statements.
\begin{enumerate}[(i)]
 \item If $D$ is a $\T$-convex set, then $\W_\T(D)$ contains all constant functions and all additive functions. Furthermore, it is closed with respect to the pointwise addition and multiplication by nonnegative scalars.
 \item If $D$ is a $\T$-convex set, then $\W_\T(D)$ is closed with respect to the pointwise chain supremum, the pointwise chain infimum, and the pointwise convergence.
\item If an endomorphism $A\in\E(X)$ commutes with any member of $\T$ and $f\in\W_\T(D)$, then $f\circ A\in\W_\T(A^{-1}(D))$.
\end{enumerate}}

\begin{proof}

The proof of (i) is obvious.

Proof of (ii). 
To prove the  $\T$-Wright convexity of the pointwise supremum of a nondecreasing family $\{f_\alpha \mid \alpha \in I\}$ of $\T$-Wright convex functions defined on $D$ fix $T \in \T$, $x, y \in D$ and $\varepsilon >0$. Let  $f:D\to[-\infty,+\infty[\,$ be given by $f=\sup \{f_\alpha \mid \alpha \in I\}$. To show that $f$ is $T$-Wright convex, observe that there exist some $\alpha_1, \alpha_2\in I$ such that 
\Eq{*}{
f(T(x)+(I-T)(y) ) &\leq f_{\alpha_1}(T(x)+(I-T)(y)) + \varepsilon \qquad \mbox{and} \\ f((I-T)(x)+T(y)) &\leq f_{\alpha_2}((I-T)(x)+T(y)) + \varepsilon.
}
Since the family $\{f_\alpha \mid \alpha \in I\}$ forms a chain, there exists $\alpha_0\in \{\alpha_1,\alpha_2\}$ such that $f_{\alpha_0}=\max(f_{\alpha_1}, f_{\alpha_2})$. Then, by the $T$-Wright convexity of $f_{\alpha_0}$, we have
\Eq{*}{
 f(T(x)+(I-T)(y)) + f((I-T)(x)+T(y)) &\leq f_{\alpha_1}(T(x)+(I-T)(y)) + f_{\alpha_2}((I-T)(x)+T(y)) + 2 \varepsilon \\&\leq
f_{\alpha_0}(T(x)+(I-T)(y)) + f_{\alpha_0}((I-T)(x)+T(y)) + 2 \varepsilon \\&\leq
 f_{\alpha_0}(x) + f_{\alpha_0}(y)+ 2\varepsilon \leq f(x) + f(y)+ 2\varepsilon.
}
Since $\varepsilon >0$ is arbitrary small, the $T$-Wright convexity of $f$ follows.

Similarly we will establish the $\T$-Wright convexity of the pointwise infimum of a chain $\{f_\alpha \mid \alpha \in I\}$ of $\T$-Wright convex functions defined on $D$. To do this, fix $T \in \T$, $x, y \in D$ and $\varepsilon >0$. Let  $f:D\to[-\infty,+\infty[\,$ be given by $f=\inf \{f_\alpha \mid \alpha \in I\}$. Then there exist some $\alpha_1, \alpha_2\in I$ such that 
$f_{\alpha_1}(x) \leq f(x) + \varepsilon$ and $f_{\alpha_2}(y) \leq f(y) + \varepsilon$. Since the family $\{f_\alpha \mid \alpha \in I\}$ forms a chain, there exists $\alpha_0\in \{\alpha_1,\alpha_2\}$ such that $f_{\alpha_0}=\min(f_{\alpha_1}, f_{\alpha_2})$ and then we have
\Eq{*}{
 f(T(x)+(I-T)(y))  + f((I-T)(x)+T(y))
 &\leq f_{\alpha_0}(T(x)+(I-T)(y)) + f_{\alpha_0}((I-T)(x)+T(y)) \\&\leq f_{\alpha_0}(x) + f_{\alpha_0}(y) 
 \leq f_{\alpha_1}(x) + f_{\alpha_2}(y) \leq f(x) + f(y) + 2\varepsilon.
}
Since $\varepsilon >0$ is arbitrary, then $f$ is $T$-Wright convex.

To show the $\T$-Wright convexity of the pointwise limit $(f_n)$ of  $\T$-Wright convex functions defined on $D$, fix $T \in \T$ and $x, y \in D$. We have 
\Eq{*}{
 f(T(x)+(I-T)(y)) + f((I-T)(x)+T(y)) &= \lim_{n \to +\infty}  f_n(T(x)+(I-T)(y)) + f_n((I-T)(x)+T(y)) \\&\leq \lim_{n \to +\infty} f_n(x)+f_n(y) = f(x)+ f(y).
}
The resulted inequality proves that $f$ is $T$-Wright convex.

Proof of (iii). To verify the $\T$-Wright convexity of the function $f\circ A$, let $x,y\in A^{-1}(D)$ and let $T\in\T$ be fixed. Then $A(x),A(y) \in D$, hence the $T$-Wright convexity of $f$ yields
\Eq{*}{
  (f\circ A)\big(T(x)+(I-T)(y)\big)&+(f\circ A)\big((I-T)(x)+T(y)\big)\\
  &=f\big(T(A(x))+(I-T)(A(y))\big) + f\big((I-T)(A(x))+T(A(y))\big)\\
  &\leq f(A(x))+f(A(y))=(f\circ A)(x)+(f\circ A)(y),
}
which completes the proof of the $\T$-Wright convexity of $f\circ A$.
\end{proof}

Now, given a function $f:D\to[-\infty,+\infty[\,$, we consider the collection of endomorphisms $T$ of $X$ that make $f$ a $T$-Wright convex function: 
\Eq{*}{
  \TW_f:=\{T\in\E(X)\mid \mbox{$f$ is $T$-Wright convex}\}.
}
It is obvious that, for every function $f$, we have $0,I\in\TW_f$ and $0,I\in\TW^d_f$ (if $(X,+,d)$ is a metric Abelian group). The next result shows some structural properties of $\TW_f$ and $\TW^d_f$.

\Thm{P1w}{Let $D\subseteq X$ and let $f:D\to[-\infty,+\infty[\,$ be an arbitrary function. Then $\TW_f$ is closed with respect to the mappings in \eq{maps}. If $(X,+,d)$ is a metric Abelian group, then $\TW^d_f$ is also closed with respect to the mappings in \eq{maps}.}

\begin{proof} The invariance of $\TW_f$ with respect to the map $T\mapsto I-T$ is an obvious consequence of the definition.

Let $T,S\in\TW_f$. Then, for all $x,y\in D$, we have
\Eq{*}{
  f(x)+f(y)
  &\geq f(S(x)+(I-S)(y))+f((I-S)(x)+S(y))\\
  &\geq f\big(T(S(x)+(I-S)(y))+(I-T)((I-S)(x)+S(y))\big)\\
  &\qquad +f\big((I-T)(S(x)+(I-S)(y))+T((I-S)(x)+S(y))\big)\\
  &=f\big((T\circ S+(I-T)\circ(I-S))(x)+(T\circ (I-S)+(I-T)\circ S)(y)\big)\\
  &\qquad +f\big((T\circ (I-S)+(I-T)\circ S)(x)+(T\circ S+(I-T)\circ(I-S))(y)\big).
}
This means that $f$ is $(T\circ S+(I-T)\circ(I-S))$-Wright convex, which was to be proved.

The proof of the second assertion is completely analogous. 
\end{proof}

The next statement is a generalization of the third assertion of Theorem 1  of the paper \cite{MakNikPal91}, which was one of the main results therein.
Our approach extensively uses \cor{nkc2} which is based on \thm{RCT}, our  generalization of the R{\aa}dstr\"om Cancellation Theorem.

\Thm{P1w+}{Let $(X,+,d)$ be a metric Abelian group, let $n_0\in\N$ such that $\mu_d(n_0)>1$, let $D$ be a closed bounded $n_0$-convex set and let $f:D\to[-\infty,+\infty[\,$ be an arbitrary function. If $n,k\in\N$, $T\in\TW^d_f$ and $S:=n\Cdot T+k\Cdot(I-T)$ is invertible with $S^{-1}\in\E^d(X)$, then $S^{-1}\circ(n\Cdot T)\in\TW^d_f$.}

\begin{proof} Let $n,k\in\N$, $T\in\TW^d_f$ be such that $S:=n\Cdot T+k\Cdot (I-T)$ is invertible with $S^{-1}\in\E^d(X)$. To prove the Wright-convexity with respect to the linear map $S^{-1}\circ(n\Cdot T)$, let $x,y\in D$ be fixed. By \cor{nkc2}, for all $(i,j)\in\{0,\dots,n\}\times\{0,\dots,k\}$, we have that $D$ is $S^{-1}\circ(i\Cdot T+j\Cdot(I-T))$-convex.
Therefore, for all $(i,j)\in\{0,\dots,n\}\times\{0,\dots,k\}$, the element $u_{i,j}$ defined by  
\Eq{*}{
  u_{i,j}:=S^{-1}\circ\big((n-i)\Cdot T+(k-j)\Cdot(I-T)\big)(x)
  +S^{-1}\circ\big(i\Cdot T+j\Cdot(I-T)\big)(y)
}
belongs to $D$. On the other hand, one can easily check that, for $(i,j)\in\{0,\dots,n-1\}\times\{0,\dots,k-1\}$, 
\Eq{*}{
  u_{i,j+1}=T(u_{i,j})+(I-T)(u_{i+1,j+1}) \qquad\mbox{and}\qquad
  u_{i+1,j}=(I-T)(u_{i,j})+T(u_{i+1,j+1}).
}
Therefore, the $T$-Wright convexity of $f$ implies that
\Eq{*}{
  f(u_{i,j+1})+f(u_{i+1,j})\leq f(u_{i,j})+f(u_{i+1,j+1})
}
for $(i,j)\in\{0,\dots,n-1\}\times\{0,\dots,k-1\}$. Adding up these inequalities side by side with respect to $i\in\{0,\dots,n-1\}$, we get
\Eq{*}{
  f(u_{0,j+1})+f(u_{n,j})\leq f(u_{0,j})+f(u_{n,j+1}).
}
Now adding up the inequalities side by side with respect to $j\in\{0,\dots,k-1\}$, we arrive at the inequality
\Eq{WW}{
  f(u_{0,k})+f(u_{n,0})\leq f(u_{0,0})+f(u_{n,k}).
}
Observe that $u_{0,0}=x$, $u_{n,k}=y$, and 
\Eq{*}{
  u_{0,k}=S^{-1}\circ(n\Cdot T)(x)+S^{-1}\circ(k\Cdot(I-T))(y)
  &=S^{-1}\circ(n\Cdot T)(x)+(I-S^{-1}\circ(n\Cdot T))(y),\\
  u_{n,0}=S^{-1}\circ(k\Cdot(I-T))(x)+S^{-1}\circ(n\Cdot T)(y)
  &=(I-S^{-1}\circ(n\Cdot T))(x)+S^{-1}\circ(n\Cdot T)(y).
}
Therefore, inequality \eq{WW} shows that $f$ is $S^{-1}\circ(n\Cdot T)$-Wright convex, which was to be proved.
\end{proof}

In the following statement, we provide conditions for the invertibility of the map $S=n\Cdot T+k\Cdot(I-T)$.

\Thm{2w}{Let $(X,+,d)$ be a complete metric Abelian group with $\mu_d(2)>1$, let $n,k\in\N$ such that $n+k\in\div(X)$ and $\mu_d(n+k)>0$, let $D$ be a closed bounded $2$-convex set, and let $f:D\to[-\infty,+\infty[\,$ be an arbitrary function. If $T\in\TW^d_f$ satisfies
\Eq{nk}{
  |n-k|\rho_d(2\Cdot T-I)<\mu_d(n+k),
}
then $X$ is $2$-divisible, $S:=n\Cdot T+k\Cdot(I-T)$ is invertible with $S^{-1}\in\E^d(X)$ and $S^{-1}\circ(n\Cdot T)\in\TW^d_f$.}

\begin{proof} For the proof of this statement, in view of \thm{P1w+}, 
it suffices to show that inequality \eq{nk} implies the invertibility of $S$ with a $d$-bounded inverse. We will prove this by using \thm{NIT}.

First observe that
\Eq{*}{
  \frac{2}{n+k}\Cdot S-I=\frac{n-k}{n+k}\Cdot(2\Cdot T-I).
}
Therefore, by the subadditivity of the $d$-norm and the submultiplicativity of $\mu_d$, we obtain
\Eq{*}{
  \Big\|\Big(\frac{2}{n+k}\Cdot S-I\Big)^m\Big\|^*_d
  =\bigg\|\frac{(n-k)^m}{(n+k)^m}\Cdot(2\Cdot T-I)^m\bigg\|^*_d
  \leq\frac{|n-k|^m}{\mu_d(n+k)^m}\|(2\Cdot T-I)^m\|^*_d.
}
Now, taking the $m$-th root side by side, then computing the upper limit as $m\to\infty$, finally using \eq{nk}, we arrive at the inequality
\Eq{*}{
  \rho_d\Big(I-\frac{2}{n+k}\Cdot S\Big)
  =\rho_d\Big(\frac{2}{n+k}\Cdot S-I\Big)
  \leq\frac{|n-k|}{\mu_d(n+k)}\rho_d(2\Cdot T-I)<1.
}
Therefore, \thm{NIT} applied for the endomorphism $I-\frac{2}{n+k}\Cdot S$ yields that $I-(I-\frac{2}{n+k}\Cdot S)=\frac{2}{n+k}\Cdot S$ is an invertible endomorphism with a bounded inverse. Thus $\pi_2$ must be a surjection and hence $2\in\div(X)$. Consequently, $\frac{1}{n+k}\Cdot S=\big(\frac{1}{n+k}\Cdot I\big)\circ S$ is also an invertible endomorphism with a bounded inverse. Therefore, $S^{-1}\in\E^d(X)$.
\end{proof}

The following result is an easy consequence of \thm{2w} and it is still more general than the third assertion of Theorem 1 of the paper \cite{MakNikPal91}.

\Cor{2w}{Let $(X,+)$ be the additive group of a Banach space, let $t\in[0,1]\cap\mathbb{Q}$, let $D$ be a closed bounded convex set, and let $f:D\to[-\infty,+\infty[\,$ be an arbitrary function. If $T\in\TW^d_f$ satisfies
\Eq{t}{
  |2t-1|\rho_d(2\Cdot T-I)<1,
}
then $S:=t\Cdot T+(1-t)\Cdot(I-T)$ is invertible with $S^{-1}\in\E^d(X)$ and $S^{-1}\circ(t\Cdot T)\in\TW^d_f$.}

\begin{proof} If $(X,+)$ is the additive group of a Banach space, then $\div(X)=\N$, $\mu_d(n)=n$ for all $n\in\N$ and convex sets are $2$-convex.
 
If $t=1$, then \eq{t} and \thm{NIT} imply that $I-(I-2\Cdot T)=2\Cdot T$ is invertible with a bounded inverse, hence $S=T$ is invertible with $S\in\E^d(X)$. The inclusion $S^{-1}\circ(t\Cdot T)=I\in\TW^d_f$ is trivial.
The case $t=0$ can analogously be seen.

In the rest of the proof, we may assume that $t\in\,]0,1[\,\cap\mathbb{Q}$. Then there exist $n,k\in\N$ such that $t=\frac{n}{n+k}$. Then inequality \eq{t} becomes \eq{nk}, hence, by \thm{2w}, we get that $(n+k)\Cdot S=n\Cdot T+k\Cdot(I-T)$ is invertible with a bounded inverse and 
\Eq{*}{
 S^{-1}\circ(t\Cdot T)
 =S^{-1}\circ\big(\tfrac{n}{n+k}\Cdot T\big)
 =((n+k)\Cdot S)^{-1}\circ (n\Cdot T)\in \TW^d_f. 
}
Therefore the proof has been completed.
\end{proof}

For an endomorphism $T\in\E(X)$, we say that a function $a:D\to[-\infty,+\infty[\,$ is \emph{$T$-Wright affine} if $D$ is a $T$-convex subset of $X$ and, for all $x,y\in D$, 
\Eq{Ta}{
  a(T(x)+(I-T)(y))+a((I-T)(x)+T(y))=a(x)+a(y).
}
If $\T\subseteq\E(X)$, then a function $a:D\to[-\infty,+\infty[\,$ is called \emph{$\T$-Wright affine} if, for all $T\in\T$, it is $T$-Wright affine. Obviously, $\T$-Wright affine functions are also $\T$-Wright convex, therefore, most of the results established above remains valid for this subclass. In what follows, we describe the solution of the above functional equation among real valued functions that are defined on the entire set $X$. Our result partially generalizes \cite[Theorem 1]{Olb17b} (cf.\ also \cite{Olb07}).

\Thm{TWA}{Let $\T\subseteq\E(X)$ and assume that there exists $T_0\in\T$  such that $T_0(X)=(I-T_0)(X)$. Then a function $a:X\to\R$ is $\T$-affine if and only if there exist a constant $c$, an additive function $A:X\to\R$ and a symmetric biadditive function $B:X\times X\to\R$ such that 
\Eq{BT}{
B(T(u),(I-T)(u))=0 \qquad(u\in X,\,T\in\T)
}
and 
\Eq{ABc}{
  a(x)=B(x,x)+A(x)+c \qquad(x\in X).
}}

\begin{proof} Assume first that $a:X\to\R$ is $\T$-Wright affine. Then, using that $a$ is real valued, by \eq{Ta} we have that
\Eq{*}{
  a(x)+a(y)-a(T_0(x)+(I-T_0)(y))-a((I-T_0)(x)+T_0(y))=0 \qquad(x,y\in X).
}
In view of the condition $T_0(X)=(I-T_0)(X)$, this functional equation
possesses the assumptions of linear functional equations dealt with by Sz\'ekelyhidi \cite[Theorem 3.6]{Sze82d}, therefore $a$ can be represented in the form \eq{ABc} for some constant constant $c$, additive function $A:X\to\R$ and symmetric biadditive function $B:X\times X\to\R$. Having this form of $a$, by the additivity of $A$, the biadditivity and symmetry of $B$, we can obtain
\Eq{*}{
  &a(x)+a(y)-a(T(x)+(I-T)(y))-a((I-T)(x)+T(y))\\
  &=B(x,x)+B(y,y)-B(T(x)\!+\!(I\!-\!T)(y),T(x)\!+\!(I\!-\!T)(y))-B((I\!-\!T)(x)\!+\!T(y),(I\!-\!T)(x)\!+\!T(y))\\
  &=2B(T(x-y),(I-T)(x-y)).
}
Therefore, the function $a$ satisfies the functional equation \eq{Ta} if and only if it has the representation \eq{ABc} and $B(T(x-y),(I-T)(x-y))=0$ holds for all $x,y\in X$, that is, if condition \eq{BT} is valid.
\end{proof}

In the paper \cite{DarLajLovMakPal07} the functional equation \eq{Ta} was considered under the assumption that $T$ is given as a multiplication by $t\in[0,1]$ and the characterization of those numbers $t$ was obtained for which there exists a nontrivial biadditive function $B$ satisfying \eq{BT}.

\section{$(T,t)$-convex and $(T,t)$-affine functions}

Assume that $(X,+)$ is an Abelian group.
For an endomorphism $T\in\E(X)$ and $t\in[0,1]$, we say that a function $f:D\to[-\infty,+\infty[\,$ is \emph{$(T,t)$-convex} if $D$ is a $T$-convex subset of $X$ and, for all $x,y\in D$, 
\Eq{Tt}{
  f(T(x)+(I-T)(y))\leq tf(x)+(1-t)f(y).
}
Here and in the rest of the paper, we use the usual convention $0\cdot(-\infty)=0$. If $\CR\subseteq\E(X)\times[0,1]$, then $f$ is called \emph{$\CR$-convex} if, for all $(T,t)\in\CR$, it is $(T,t)$-convex. In particular, if $\T\subseteq\E(X)$ and $\tau:\T\to[0,1]$, then $f$ is called \emph{$(\T,\tau)$-convex} if, for all $T\in\T$, it is $(T,\tau(T))$-convex, that is, if $f$ is $\CR$-convex where $\CR:=\{(T,\tau(T))\mid T\in\T\}$. The class of $\CR$-convex, in particular, $(\T,\tau)$-convex functions defined on $D$ are denoted by $\C_\CR(D)$ and $\C_{\T,\tau}(D)$, respectively. If $(X,+)$ is the additive group of a linear space and $f$ is $(t\Cdot I,t)$-convex for some $t\in[0,1]$, then we say that $f$ is \emph{$t$-convex}. If $f$ is $t$-convex for all $t\in[0,1]$, then it is called a \emph{convex} function (in the standard sense). If $(X,+)$ is a uniquely 2-divisible group, then $(\frac12\Cdot I,\frac12)$-convex functions are exactly the Jensen convex ones. 

One can observe that the function $f\equiv-\infty$ is trivially $(T,t)$-convex for arbitrary $(T,t)\in \E(X)\times[0,1]$. On the other hand, it is possible that a $(T,t)$-convex function can take both finite and infinite values. To exclude this possibility, the following lemma will be useful. In what follows, for a $T$-convex set $D$, we say that an element $p\in D$ is $T$-internal with respect to $D$, if $D$ has no proper subset $E$ which contains $p$ and, for all $x,y\in D$ with $T(x)+(I-T)(y)\in E$ implies $x,y\in E$.

\Lem{inf}{Let $(T,t)\in\E(X)\times\,]0,1[\,$ and let $D\subseteq X$ be a $T$-convex set. If $f:D\to[-\infty,+\infty[\,$ is $(T,t)$-convex, then either $f(x)=-\infty$ for all $T$-interior points $x\in D$ or $f(x)\in\R$ for all $x\in D$. }

\begin{proof}
Assume that $f$ is a $(T,t)$-convex function which is finite at some $T$-interior point $p$. Define the sequence of sets $(D_n)$ by the recursion
\Eq{Dn}{
  D_0:=\{p\}, \qquad 
  D_n:=\bigcup\{\{x,y\}\mid x,y\in D,\, T(x)+(I-T)(y)\in D_{n-1}\} \qquad(n\in\N).
}
Observe that $p\in D_1$, which implies $D_0\subseteq D_1$, and hence $(D_n)$ is an increasing sequence of sets. Let $E:=\bigcup_{n=0}^\infty D_n$. Then $p\in E$, and taking the union of both sides in \eq{Dn}, it follows that
\Eq{*}{
  E=\bigcup\{\{x,y\}\mid x,y\in D,\,  T(x)+(I-T)(y)\in E\}.
}
Therefore, $E$ is a subset of $D$ which contains $p$ and, for all $x,y\in D$ with $T(x)+(I-T)(y)\in E$ implies $x,y\in E$. By the $T$-internality of $p$, it follows that $E=D$, that is, $D=\bigcup_{n=0}^\infty D_n$.

In the rest of the proof, we show that $f$ is finite-valued on $D_n$ for all $n\geq0$. This is obvious for $n=0$ by the choice of $p$. Now assume that $f$ is finite-valued on $D_{n-1}$ for some $n\in\N$. Let $z\in D_n$. Then there exist $x,y\in D$ such that $z\in\{x,y\}$ and $T(x)+(I-T)(y)\in D_{n-1}$.
Then $f(T(x)+(I-T)(y))>-\infty$ and, by the $(T,t)$-convexity of $f$,
\eq{Tt} holds. The left hand side being finite, the condition $t(1-t)>0$ implies that $f(x)$ as well as $f(y)$ are also finite, which yields that $f$ is finite at $z$. This completes the induction and finally shows that $f$ is finite valued on $D$.
\end{proof}

\Lem{IX}{Let $T\in\E(X)$ such that either $T(X)\subseteq (I-T)(X)$ or $(I-T)(X)\subseteq T(X)$ holds. Then every element of $X$ is $T$-internal with respect to $X$.}

\begin{proof} Let $p\in X$ and let $E$ be a set which contains $p$ and, for all $x,y\in D$ with $T(x)+(I-T)(y)\in E$ implies $x,y\in E$. Then
\Eq{*}{
  D_1=\bigcup\{\{x,y\}\mid x,y\in X,\, T(x)+(I-T)(y)=p\}\subseteq E
}
We show that $D_1=X$, which shows that $E$ cannot be proper and hence $p$ must be $T$-internal with respect to $X$.

Assume that $T(X)\subseteq (I-T)(X)$ holds and let $x\in X$ be arbitrary. Then $T(x-p)\in (I-T)(X)$, therefore, there exists $y\in X$ such that $T(x-p)=(I-T)(y)$. Hence, $T(x)+(I-T)(y)=T(p)+(I-T)(p)=p$, which shows that $x\in D_1$. Thus, in this case we have obtained that $D_1=X$. In the other case, the argument is completely analogous.
\end{proof}

The epigraph of an arbitrary function $f:D\to[-\infty,+\infty[\,$ is defined by
\Eq{*}{
  \epi(f):=\{(x,u)\in D\times\R\mid f(x)\leq u\}.
}
For $T\in\E(X)$ and $t\in\R$, the endomorphism $(T,t)\in\E(X\times\R)$ is defined as
\Eq{*}{
  (T,t)(x,u):=(T(x),t\Cdot u)\qquad((x,u)\in X\times\R).
}
Therefore, any relation $\CR\subseteq\E(X)\times\R$ can be viewed as a subset of $\E(X\times\R)$ as well. For $\CR\subseteq\E(X)\times\R$, we introduce the domain and codomain of $\CR$ as follows
\Eq{*}{
  \dom(\CR):=\{T\in\E(X)\mid \exists t\in\R:(T,t)\in\CR\},\qquad
  \codom(\CR):=\{t\in\R\mid \exists T\in\E(X):(T,t)\in\CR\}.
}

The following characterization of $\CR$-convexity of functions is important.

\Thm{ChC}{Let $\CR\subseteq\E(X)\times[0,1]$ be a nonempty subset. Then a function $f:D\to[-\infty,+\infty[\,$ is $\CR$-convex if and only if, its epigraph $\epi(f)$ is an $\CR$-convex subset of $X\times\R$.}

\begin{proof}
To prove that the epigraph of an $\CR$-convex function $f$ is an $\CR$-convex set, let $(T,t) \in \CR$. Then $D$ is a $T$-convex subset of $X$. 
Fix some $p, q \in \epi(f)$. There exist $x, y \in D$ and $u, v \in \R$ such that $p=(x,u)$, $q=(y,v)$ and $f(x)\leq u$, $f(y)\leq v$. From the $T$-convexity of $D$ and from the $(T,t)$-convexity of $f$, we get
\Eq{*}{
f\big(  T(x) + (I-T)(y)  \big) \leq tf(x) + (1-t)f(y) \leq tu + (1-t)v.
}
This inequality gives us 
\Eq{*}{
(T,t)(x,u)+(I-T,1-t)(y,v)= \big(T(x) + (I- T)(y), tu +(1- t)v \big)  \in \epi(f),
}
which shows that $\epi(f)$ is $(T,t)$-convex.

To prove the converse implication, let $(T,t) \in \CR$ and $x, y \in D$. We have $(x,f(x)), (y,f(y)) \in \epi(f)$. Therefore, by the assumed $(T,t)$-convexity of $\epi(f)$, 
\Eq{*}{
\big( T(x)  + (I-T)(y), tf(x) + (1-t)f(y)  \big) 
= (T,t)(x,f(x))+(I-T,1-t)(y,f(y)) \in \epi(f),
}
which yields the $T$-convexity of $D$ and the $(T,t)$-convexity of $f$, and completes the proof.
\end{proof}

\Thm{TC}{Let $\CR\subseteq\E(X)\times[0,1]$ be a nonempty subset. Then we have the following statements.
\begin{enumerate}[(i)]
 \item $\C_{\CR}(D)$ is closed with respect to the pointwise addition and multiplication by nonnegative scalars.
 \item $\C_{\CR}(D)$ is closed with respect to the pointwise supremum, the pointwise chain infimum, and the pointwise convergence.
\item If $f:D\to[-\infty,+\infty[\,$ and $g:E\to[-\infty,+\infty[\,$ are $\CR$-convex functions, then the function $f*g:(D+E)\to[-\infty,+\infty[\,$ defined by
\Eq{*}{
  (f*g)(x):=\inf\{f(u)+g(v)\mid u\in D,\,v\in E,\,u+v=x\}
}
is $\CR$-convex on $D+E$.
\item If $A\in\E(X)$ commutes with any member of the domain of $\CR$ and $f\in\C_{\CR}(D)$, then $f\circ A\in\C_{\CR}(A^{-1}(D))$ and $f\circ A^{-1}\in\C_{\CR}(A(D))$, where $f\circ A^{-1}:A(D)\to[-\infty,+\infty[\,$ is defined by \eq{fAi}.
\end{enumerate}}

\begin{proof}
The proof of (i) is obvious.

Proof of (ii). 
To verify the  $\CR$-convexity of the pointwise supremum of a
family $\{f_\alpha \mid \alpha \in I\}$ of $\CR$-convex functions defined on $D$, fix $(T,t) \in \CR$, $x, y \in D$ and $\varepsilon >0$. Let  $f:D\to[-\infty,+\infty[\,$ be given by $f=\sup \{f_\alpha \mid \alpha \in I\}$. To show that $f$ is $(T,t)$-convex observe that there exists some $\alpha_0\in I$ such that 
\Eq{*}{
f(T(x)+(I-T)(y) ) &\leq f_{\alpha_0}(T(x)+(I-T)(y)) + \varepsilon \\&\leq t f_{\alpha_0}(x) + (1-t)f_{\alpha_0}(y) + \varepsilon \leq t f(x) + (1-t)f(y) + \varepsilon.
}
Since $\varepsilon >0$ is arbitrary small, therefore $f$ is $(T,t)$-convex.

To justify the $\CR$-convexity of the pointwise infimum of a chain $\{f_\alpha \mid \alpha \in I\}$ of $\CR$-convex functions defined on $D$, let $(T,t) \in \CR$, $x, y \in D$ and $\varepsilon >0$. Let $f$ be given by $f=\inf \{f_\alpha \mid \alpha \in I\}$. Then there exist $\alpha_1, \alpha_2\in I$ such that $f_{\alpha_1}(x) \leq f(x) + \varepsilon$ and $f_{\alpha_2}(y) \leq f(y) + \varepsilon$. By the chain property, there exists $\alpha_0\in \{\alpha_1,\alpha_2\}$ such that $f_{\alpha_0}=\min(f_{\alpha_1}, f_{\alpha_2})$. Then we have
\Eq{*}{
 f(T(x)+(I-T)(y)) &\leq f_{\alpha_0}(T(x)+(I-T)(y))  
 \leq tf_{\alpha_0}(x)+ (1-t)f_{\alpha_0}(y)
 \leq tf_{\alpha_1}(x)+ (1-t)f_{\alpha_2}(y) \\
 &\leq t(f(x) + \varepsilon)+ (1-t)(f(y)+ \varepsilon ) 
 = tf(x) + (1-t)f(y)  + \varepsilon. 
}
Upon taking the limit $\varepsilon \to0$, we obtain that $f$ is $(T,t)$-convex.

To show the $\CR$-convexity of the pointwise limit $(f_n)$ of $\CR$-convex functions defined on $D$, let $(T,t) \in \CR$ and $x, y \in D$. We have 
\Eq{*}{
 f(T(x)+(I-T)(y)) &= \lim_{n \to +\infty}  f_n(T(x)+(I-T)(y))\\ 
 &\leq \lim_{n \to +\infty} tf_n(x)+ (1-t)f_n(y) =  tf(x)+ (1-t)f(y) .
}
This inequality proves that $f$ is $(T,t)$-convex.

Proof of (iii). For the $\CR$-convexity of the function $f*g$, let $(T,t) \in \CR$ and $x,y\in D+E$. We need to prove that 
\Eq{f**g}{
  (f*g)(T(x)+(I-T)(y))\leq t(f*g)(x)+(1-t)(f*g)(y).
}
Let $c\in \,](f*g)(x),+\infty[\,$ and $d\in \,](f*g)(y),+\infty[\,$ be arbitrary. Then, by the definition of $f*g$, there exist $u,v\in D$ such that
\Eq{*}{
  f(u)+g(x-u)<c \qquad\mbox{and}\qquad f(v)+g(y-v)<d.
}
Then, by the $(T,t)$-convexity of $f$ and $g$, we obtain
\Eq{*}{
  (f*g)(T(x)+(I-T)(y))
  &\leq f(T(u)+(I-T)(v))+g(T(x-u)+(I-T)(y-v))\\
  &\leq (tf(u)+(1-t)f(v))+(tg(x-u)+(1-t)g(y-v)) \\
  &= t(f(u)+g(x-u))+(1-t)(f(v)+g(y-v))\\
  &<tc+(1-t)d.
}
Upon taking the limits $c\searrow (f*g)(x)$ and $d\searrow (f*g)(y)$, the inequality \eq{f**g} follows.
 
Proof of (iv). To verify the $\CR$-convexity of the function $f\circ A$, let $x,y\in A^{-1}(D)$ and $(T,t)\in\CR$. Then $A(x),A(y) \in D$, hence the $(T,t)$-convexity of $f$ yields
\Eq{*}{
  (f\circ A)(T(x)+(I-T)(y))
  &=f\big(T(A(x))+(I-T)(A(y))\big)\\
  &\leq tf(A(x))+(1-t)f(A(y))
  =t(f\circ A)(x)+(1-t)(f\circ A)(y),
}
which proves the $(T,t)$-convexity of $f\circ A$.

Finally, we show the $\CR$-convexity of the function $f\circ A^{-1}$. For this proof, let $x,y\in A(D)$ and $(T,t)\in\CR$. For the proof of the inequality 
\Eq{fA}{
  (f\circ A^{-1})(T(x)+(I-T)(y))\leq t(f\circ A^{-1})(x)+(1-t)(f\circ A^{-1})(y)
}
choose $c\in \,](f\circ A^{-1})(x),+\infty[\,$ and $d\in \,](f\circ A^{-1})(y),+\infty[\,$ arbitrarily. Then, there exist $u,v\in D$ such that $A(u)=x$, $A(v)=y$ and $f(u)<c$, $f(v)<d$. Using the equality $A(T(u)+(I-T)(v))=T(x)+(I-T)(y)$ and the $(T,t)$-convexity of $f$, we get
\Eq{*}{
  (f\circ A^{-1})(T(x)+(I-T)(y))
  &\leq f(T(u)+(I-T)(v))\\
  &\leq tf(u)+(1-t)f(v)
  <tc+(1-t)d.
}
Upon taking the limits $c\searrow (f\circ A^{-1})(x)$ and $d\searrow (f\circ A^{-1})(y)$, the inequality \eq{fA} follows.
\end{proof}

Now, given a function $f:D\to[-\infty,+\infty[\,$, we consider the set of those pairs $(T,t)\in\E(X)\times[0,1]$ that make $f$ a $(T,t)$-convex function: 
\Eq{*}{
  \CR_f:=\{(T,t)\in\E(X)\times[0,1]\mid \mbox{$f$ is $(T,t)$-convex}\}.
}
It is easy to see that, for all $T\in\E(X)$, the set $\CR_f(T)$ is a closed (possibly empty) subinterval of $[0,1]$. Obviously, $0\in\CR_f(0)$ and $1\in\CR_f(I)$ for every function $f:D\to[-\infty,+\infty[\,$. If $f$ is nonconstant, then the equalities $\CR_f(0)=\{0\}$ and $\CR_f(I)=\{1\}$ can easily be seen. 

It is obvious that, for every function $f$, we have $0,I\in\dom(\CR_f)$ and $0,I\in(\dom(\CR_f))^d$ (if $(X,+,d)$ is a metric Abelian group). The next result shows some structural properties of $\dom(\CR_f)$ and $(\dom(\CR_f))^d$.

\Thm{P1c}{Let $D\subseteq X$ and let $f:D\to[-\infty,+\infty[\,$. Then $\CR_f$ is a $\widetilde{\CR}_f$-convex subset of $\E(X)\times[0,1]$. If $(X,+,d)$ is a metric Abelian group, then $\CR_f^d$ is a $\widetilde{\CR}_f^d$-convex subset of $\E^d(X)\times[0,1]$. In particular, these sets are closed with respect to (componentwise) multiplication.}

\begin{proof} Let $(T,t),(T_1,t_1),(T_2,t_2)\in\CR_f$. Set 
\Eq{*}{
  S:=T\circ T_1+(I-T)\circ T_2 \qquad\mbox{and}\qquad 
  s:=tt_1+(1-t)t_2.
}
Then, $f$ is $(T_1,t_1)$-, $(T_2,t_2)$- and $(T,t)$-convex, therefore, for all $x,y\in D$, we have
\Eq{*}{
  f(T_1(x)+(I-T_1)(y))&\leq t_1f(x)+(1-t_1)f(y), \\
  f(T_2(x)+(I-T_2)(y))&\leq t_2f(x)+(1-t_2)f(y).
}
Consequently,
\Eq{*}{
  f(S(x)+(I-S)(y))
  &=f\big(T(T_1(x)+(I-T_1)(y))+(I-T)(T_2(x)+(I-T_2)(y))\big)\\
  &\leq tf\big(T_1(x)+(I-T_1)(y)\big)
     +(1-t)f\big(T_2(x)+(I-T_2)(y)\big) \\
  &\leq t\big(t_1 f(x)+(1-t_1)f(y)\big)
     +(1-t)\big(t_2f(x)+(1-t_2)f(y)\big) \\
  &\leq sf(x)+(1-s)f(y).
}
This means that $f$ is $(S,s)$-convex, which was to be proved.

The proof of the second assertion is completely analogous. The last assertion easily follows from the first two by taking $(T_2,t_2)=(0,0)$ in the above proof.
\end{proof}

If $(T,t)\in\E(X)\times[0,1]$, then we say that a function $a:D\to[-\infty,+\infty[\,$ is \emph{$(T,t)$-affine} if $D$ is a $T$-convex subset of $X$ and, for all $x,y\in D$, 
\Eq{AF}{
  a(T(x)+(I-T)(y))=ta(x)+(1-t)a(y).
}
If $\CR\subseteq\E(X)\times[0,1]$, then a function $a:D\to[-\infty,+\infty[\,$ is called \emph{$\CR$-affine} if, for all $(T,t)\in\CR$, it is $(T,t)$-affine. The class of $\CR$-affine functions defined on $D$ is denoted by $\A_{\CR}(D)$. If $(X,+)$ is the additive group of a linear space and $f$ is $(t\Cdot I,t)$-affine for some $t\in[0,1]$, then we say that $f$ is \emph{$t$-affine}. If $f$ is $t$-affine for all $t\in[0,1]$, then it is called an \emph{affine} function. If $(X,+)$ is a uniquely 2-divisible group, then $(\frac12\Cdot I,\frac12)$-affine functions are called \emph{Jensen affine} or \emph{midpoint affine}.

The graph of a function $a:D\to[-\infty,+\infty[\,$ is defined by
\Eq{*}{
  \graph(f):=\{(x,a(x))\mid x\in D,\, -\infty<a(x)\}.
}
The following characterization of $\CR$-affinity is important.

\Thm{ChA}{Let $\CR\subseteq\E(X)\times[0,1]$ be a nonempty subset. Then a function $a:D\to\R$ is $\CR$-affine if and only if $\graph(a)$ is an $\CR$-convex subset of $X\times\R$.}

\begin{proof}
To prove that the graph of an $\CR$-affine function $a$ is an $\CR$-convex set, let $(T,t) \in \CR$. Then $D$ is a $T$-convex subset of $X$. 
Fix some $p, q \in \graph(a)$. Then there exist $x,y \in D$ such that$p=(x,a(x))$, $q=(y,a(y))$. From the $T$-convexity of $D$ and from the $(T,t)$-affinity of $a$, we have \eq{AF},
which gives us 
\Eq{*}{
(T,t)(x,a(x))+(I-T,1-t)(y,a(y))
&= \big(T(x) + (I- T)(y), t a(x) +(1- t) a(y) \big) \\
&= \big(T(x) + (I- T)(y), a(T(x) + (I- T)(y)) \big) \in \graph(a),
}
which shows that $\graph(a)$ is $(T,t)$-convex.

To prove the converse implication, let $(T,t) \in \CR$ and $x, y \in D$. We have $(x,f(x)), (y,f(y)) \in \graph(a)$. Therefore, by the assumed $(T,t)$-convexity of $\graph(a)$, 
\Eq{*}{
\big( T(x)  + (I-T)(y), ta(x) + (1-t)a(y)  \big) 
= (T,t)(x,a(x))+(I-T,1-t)(y,a(y)) \in \graph(a),
}
which yields the $T$-convexity of $D$ and the $(T,t)$-affinity of $a$, and completes the proof.
\end{proof}

\Thm{TA}{Let $\CR\subseteq\E(X)\times[0,1]$ be a nonempty subset. Then we have the following statements.
\begin{enumerate}[(i)]
 \item $\A_{\CR}(D)$ is closed with respect to the pointwise addition and multiplication by nonnegative scalars.
 \item $\C_{\CR}(D)$ is closed with respect to the pointwise convergence.
\item If $A\in\E(X)$ commutes with any member of the domain of $\CR$ and $a\in\A_{\CR}(D)$, then $a\circ A\in\A_{\CR}(A^{-1}(D))$.
\end{enumerate}}

\begin{proof}
The proof of (i) is obvious and the proofs of (ii) and (iii) are parallel to the corresponding statements of \thm{TC}, therefore, they are omitted.
\end{proof}

Now, given $a:D\to[-\infty,+\infty[\,$, we consider the set of those pairs $(T,t)\in\E(X)\times[0,1]$ such that $a$ is a $(T,t)$-affine function: 
\Eq{*}{
  \S_a:=\{(T,t)\in\E(X)\times[0,1]\mid \mbox{$a$ is $(T,t)$-affine}\}.
}
Obviously, $\S_a\subseteq\CR_a$ and, for all $T\in\E(X)$, the set $\S_a(T)$ is a closed (possibly empty) subinterval of $[0,1]$. In addition, $0\in\S_0(0)$ and $1\in\S_0(I)$ for every function $a:D\to[-\infty,+\infty[\,$. If $a$ is nonconstant, then the equalities $\S_a(0)=\{0\}$ and $\S_a(I)=\{1\}$ can easily be seen. 

The next result establishes some structural properties of $\dom(\S_a)$ and $(\dom(\S_a))^d$.

\Thm{P1A}{Let $D\subseteq X$ and let $a:D\to[-\infty,+\infty[\,$. Then $\S_a$ is a $\widetilde{\S}_a$-convex subset of $\E(X)\times[0,1]$. If $(X,+,d)$ is a metric Abelian group, then $\S_a^d$ is a $\widetilde{\S}_a^d$-convex subset of $\E^d(X)\times[0,1]$. In particular, these sets are closed with respect to (componentwise) multiplication.}

\begin{proof} Let $(T,t),(T_1,t_1),(T_2,t_2)\in\S_a$. Set 
\Eq{*}{
  S:=T\circ T_1+(I-T)\circ T_2 \qquad\mbox{and}\qquad 
  s:=tt_1+(1-t)t_2.
}
Then, $a$ is $(T_1,t_1)$-, $(T_2,t_2)$- and $(T,t)$-affine, therefore, for all $x,y\in D$, we have
\Eq{*}{
  a(T_1(x)+(I-T_1)(y))&= t_1a(x)+(1-t_1)a(y), \\
  a(T_2(x)+(I-T_2)(y))&= t_2a(x)+(1-t_2)a(y).
}
Consequently,
\Eq{*}{
  a(S(x)+(I-S)(y))
  &=a\big(T(T_1(x)+(I-T_1)(y))+(I-T)(T_2(x)+(I-T_2)(y))\big)\\
  &= ta\big(T_1(x)+(I-T_1)(y)\big)
     +(1-t)a\big(T_2(x)+(I-T_2)(y)\big) \\
  &= t\big(t_1 a(x)+(1-t_1)a(y)\big)
     +(1-t)\big(t_2a(x)+(1-t_2)a(y)\big) \\
  &= sa(x)+(1-s)a(y).
}
This means that $a$ is $(S,s)$-affine, which was to be proved.

The proof of the second assertion is analogous. The last assertion easily follows from the first two by taking $(T_2,t_2)=(0,0)$ in the above proof.
\end{proof}

\Thm{P2A}{Let $D\subseteq X$, $a:D\to\R$ and $(T,t),(S,s)\in\S_a$ with $0<s$, $t\leq s$ and $S^*:T(X)\to X$ be a homomorphism which is the right inverse of $S$ on the codomain of $T$, i.e., $S\circ S^*(u)=u$ for all $u\in T(X)$. Then $a$ is $(S^*\circ T,s^{-1}t)$- and $(S-T,s-t)$-affine.}

\begin{proof}
The proof is based on the following identity:
\Eq{*}{
 S\big(S^*\circ T(x)+(I-S^*\circ T)(y)\big)+(I-S)(y)
 = T(x)+(I-T)(y).
}
Using this, the $(S,s)$- and $(T,t)$-affinity of $a$ yield
\Eq{*}{
  sa\big(S^*\circ T(x)+(I-S^*\circ T)(y)\big)+(1-s)a(y)
  &=a\big(S\big(S^*\circ T(x)+(I-S^*\circ T)(y)\big)+(I-S)(y)\big)\\
  &=a\big(T(x)+(I-T)(y)\big)=ta(x)+(1-t)a(y).\\
}
Subtracting $(1-s)a(y)$ from both sides and then dividing the inequality so obtained by $s$, we get
\Eq{*}{
  a\big(S^*\circ T(x)+(I-S^*\circ T)(y)\big)
  =(s^{-1}t)a(x)+(1-s^{-1}t)a(y).
}
This proves that $a$ is $(S^*\circ T,s^{-1}t)$-affine.

To show the last assertion, observe that $a$ is $(I-S^*\circ T,1-s^{-1}t)$-affine and
\Eq{*}{
  (S-T,s-t)=(S\circ (I-S^*\circ T),s(1-s^{-1}t)).
}
Thus, the statement follows from \thm{P1A}.
\end{proof}

It follows from the above theorem that if $S^*:(I-T)(X)\to X$ is a homomorphism which is the right inverse of $S$ on the codomain of $I-T$ and $s+t\geq1$, then $a$ is also $(S+T-I,s+t-1)$-affine. Indeed, 
\Eq{*}{
  (S+T-I,s+t-1)=(S-(I-T),s-(1-t)).
}
Using that $(I-T,1-t)\in\S_a$, the above equality and the last assertion of \thm{P2A} imply that $a$ is also $(S+T-I,s+t-1)$-affine.

Obviously, the function $a\equiv-\infty$ is $(T,t)$-affine for all $(T,t)\in\E(X)\times[0,1]$. The following statement describes a large class of nontrivial real valued $(T,t)$-affine functions.

\Prp{AF}{Let $(T,t)\in\E(X)\times[0,1]$, $A\in\E(X)$ and $c\in\R$ be such that $A$ is $(T,t)$-homogeneous, i.e., $A\circ T=tA$ holds. Then $a:=A+c$ is $(T,t)$-affine. Furthermore, if $\CR\subseteq\E(X)\times[0,1]$ is nonempty, and $A$ is $(T,t)$-homogeneous for all pairs $(T,t)\in\CR$, then $a$ is $\CR$-affine.}

\begin{proof} Let $x,y\in X$, then
\Eq{*}{
 a(T(x)+(I-T)(y))
 &=A(T(x)+(I-T)(y))+c=A(T(x-y)+y)+c\\
 &=A(T(x-y))+A(y)+c=tA(x-y)+A(y)+c\\
 &=t(A(x)-A(y))+A(y)+c=t(A(x)+c)+(1-t)(A(y)+c)\\
 &=ta(x)+(1-t)a(y),
}
which shows that \eq{AF} holds proving that $a$ is $(T,t)$-affine.
The last statement immediately follows from the first part.
\end{proof}

The following theorem offers a characterization of $\CR$-affine functions defined on $X$. The main tool for the proof is the result of Sz\'ekelyhidi \cite{Sze82d}, which describes the general solution of linear functional equations with constant coefficients.

\Thm{AF}{Let $\CR\subseteq\E(X)\times[0,1]$ and assume that there exists $(T_0,t_0)\in\CR$ such that $t_0\in\,]0,1[\,$ and either $T_0(X)\subseteq (I-T_0)(X)$ or $(I-T_0)(X)\subseteq T_0(X)$ holds. If $a:X\to[-\infty,+\infty[\,$ is $\CR$-affine, then either $a\equiv-\infty$ or there exist an additive function $A:X\to\R$ and $c\in\R$ such that $A$ is $(T,t)$-homogeneous for all pairs $(T,t)\in\CR$ and $a=A+c$ holds.}

\begin{proof} Assume that $a$ is an $\CR$-affine function which is not identically equal to $-\infty$. First we show that $a(x)\in\R$ for all $x\in X$. There exists $p\in X$ such that $-\infty<a(p)$, i.e., $a(p)\in\R$. By \lem{IX}, it follows that $p$ is $T_0$-internal with respect to $X$, therefore, by \lem{inf}, $a$ is finite everywhere.

In the case when $T_0(X)\subseteq (I-T_0)(X)$ holds, we can rewrite the $(T_0,t_0)$-affinity of $a$ in the following form
\Eq{*}{
  a(x)+\frac{1-t_0}{t_0}a(y)-\frac{1}{t_0}a(T_0(x)+(I-T_0)(y))=0 \qquad(x,y\in X).
}
This is a particular case of the linear functional equations investigated by Sz\'ekelyhidi. Therefore, by  \cite[Theorem 3.6]{Sze82d}, it follows that $a$ is a first-degree generalized polynomial, that is, $a=A+c$, where $A\in\E(X)$ and $c\in\R$. Now the $\CR$-affine property of $a$ implies that
\Eq{*}{
  A(T(x)+(I-T)(y))=tA(x)+(1-t)A(y) \qquad(x,y\in X)
}
holds for all $(T,t)\in\CR$. Putting $y=0$ in this equality, it follows that $A$ is $(T,t)$-homogeneous.

The case when $(I-T_0)(X)\subseteq T_0(X)$ is valid, is completely analogous.
\end{proof}

The following characterization of $\CR$-convex functions is based on Rod\'e's celebrated separation theorem \cite{Rod78}. A relation $\CR\subseteq\E(X)\times[0,1]$ will be called nonsingular if $(T,0)\in\CR$ implies $T=0$ and $(T,1)\in\CR$ implies $T=I$.

\Thm{Rod}{Let $D\subseteq X$, let $\CR\subseteq\E(X)\times[0,1]$ be a nonempty nonsingular relation such that $\dom(\CR)$ forms a pairwise commuting subfamily of $\E(X)$. Then a function $f:D\to[-\infty,+\infty[\,$ is $\CR$-convex if and only if, for all $p\in D$, there exists an $\CR$-affine function $a:D\to[-\infty,+\infty[\,$ such that
\Eq{af}{
  a(p)=f(p)\qquad\mbox{and}\qquad a\leq f.
}}

\begin{proof} Assume first that, for all $p\in D$, there exists an $\CR$-affine function $a:D\to[-\infty,+\infty[\,$ satisfying \eq{af}. Then, $f$ is the pointwise supremum of $\CR$-affine and hence $\CR$-convex functions. Thus, by the second assertions of \thm{TC}, it is an $\CR$-convex function.

To prove the other implication, assume that $f$ is $\CR$-convex and, for an endomorphism $T\in\dom(\CR)$, define the binary operation $\omega_T:X^2\to X$ by
\Eq{*}{
  \omega_T(x,y):=T(x)+(I-T)(y).
}
Then, $\omega_T$ is idempotent and, by the $\dom(\CR)$-convexity of $D$, it follows that $D$ is closed with respect to the operation $\omega_T$ for all $T\in\dom(\CR)$.
The assumption that $\dom(\R)$ forms a pairwise commuting family of endomorphisms, easily implies that the set of operations $\{\omega_T\mid T\in\dom(\CR)\}$ is also commuting in the following sense:
\Eq{*}{
  m_T(m_S(x,y),m_S(u,v))=m_S(m_T(x,u),m_T(y,v))
  \qquad(x,y,u,v\in X,\,T,S\in\dom(\CR)).
}
Therefore, all basic assumptions of the theorem of Rod\'e are satisfied. The $\CR$-convexity of $f$ is now equivalent to the property
\Eq{*}{
  f(m_T(x,y))\leq tf(x)+(1-t)f(y) \qquad(x,y\in D,\, (T,t)\in\CR).
}
Let now $p\in D$ be fixed and define $g:D\to[-\infty,+\infty[\,$ by
\Eq{*}{
  g(p):=f(p)\qquad\mbox{and}\qquad g(x):=-\infty \quad(x\in D\setminus\{p\}).
}
Then, by the idempotent property of the operation $\omega_T$ and by the nonsingularity of the relation $\CR$, we can see that $g$ satisfies the inequality,
\Eq{*}{
  g(m_T(x,y))\geq tg(x)+(1-t)g(y) \qquad(x,y\in D,\, (T,t)\in\CR),
}
i.e., $g$ is $\CR$-concave. In addition, we trivially have that $g\leq f$ on $D$. Thus, by the theorem of Rod\'e, it follows that there exists a function $a:D\to[-\infty,+\infty[\,$ between $g$ and $f$ which satisfies
the equality
\Eq{*}{
  a(m_T(x,y))= ta(x)+(1-t)a(y) \qquad(x,y\in D,\, (T,t)\in\CR),
}
which means that $a$ is $\CR$-affine. The inequalities $f(p)=g(p)\leq a(p)\leq f(p)$ yield that $a(p)=f(p)$.
\end{proof}

The following result is a direct consequence of \thm{AF} and \thm{Rod}.

\Cor{Rod}{Let $\CR\subseteq\E(X)\times[0,1]$ be a nonsingular relation such that $\dom(\CR)$ forms a pairwise commuting subfamily of $\E(X)$ and assume that there exists $(T_0,t_0)\in\CR$ such that $t_0\in\,]0,1[\,$ and either $T_0(X)\subseteq (I-T_0)(X)$ or $(I-T_0)(X)\subseteq T_0(X)$ holds. Then a function $f:X\to[-\infty,+\infty[\,$ is $\CR$-convex if and only if, either $f\equiv-\infty$ or, for all $p\in X$, there exist an additive function $A:X\to\R$ and a constant $c\in\R$ such that $A$ is $(T,t)$-homogeneous for all pairs $(T,t)\in\CR$ and
\Eq{*}{
  A(p)+c=f(p)\qquad\mbox{and}\qquad A+c\leq f.
}}

\Thm{Last}{Let $(X,+,d)$ be a metric  Abelian group, let $n_0\in\N$ such that $\mu_d(n_0)>1$ and let $D$ be a closed bounded $n_0$-convex set. Let $f:D\to\R$, $n\in\N$ and $(T_1,t_1),\dots,(T_n,t_n)\in\CR_f^d$ with $t_1,\dots,t_n\in\,]0,1[\,$. Assume that the endomorphisms $T_1,\dots,T_n$ are pairwise commuting and, for $k\in\{0,\dots,n\}$, define
\Eq{S-s}{
  S_k:= T_1\circ \dots\circ T_k\circ 
  (I-T_{k+1})\circ\dots\circ(I-T_n) \qquad\mbox{and}\qquad
  s_k:= t_1\cdots t_k\cdot (1-t_{k+1})\cdots(1-t_n).
}
Assume that $S:=S_0+\dots+S_n$ is a bijection with $S^{-1}\in\E^d(X)$ and denote $s:=s_0+\dots+s_n$. Then, for all $k\in\{1,\dots,n\}$, we have $\big(S^{-1}\circ(S_k+\cdots+S_n),s^{-1}(s_k+\dots+s_n)\big)\in\CR_f^d$.}

\begin{proof} According to the last assertion of \thm{P1}, we have that $S_k\in\T_D^d$ for each $k\in\{0,\dots,n\}$. Thus, using that $D$ be a closed bounded $n_0$-convex set for some $n_0$ with $\mu_d(n_0)>1$ and applying \cor{nkc2}, we obtain $R_k:=S^{-1}\circ(S_k+\cdots+S_n)\in\T_D^d$ also for each $k\in\{0,\dots,n\}$. We also denote $r_k:=s^{-1}(s_k+\dots+s_n)$. We can see that $(R_0,r_0)=(I,1)$ and for the sake of brevity, let $(R_{n+1},r_{n+1}):=(0,0)$.

Using the commuting property of the endomorphisms, for all $i\in\{1,\dots,n\}$, we have that
\Eq{Ss}{
  T_i\circ S_{i-1}=(I-T_i)\circ S_i \qquad\mbox{and similarly}\qquad 
  t_is_{i-1}=(1-t_i)s_i.
}
Therefore,
\Eq{*}{
  T_i\circ (S_{i-1}+\dots+S_n)&+(I-T_i)\circ (S_{i+1}+\dots+S_n)\\
  &=T_i\circ S_{i-1}+T_i\circ S_i+S_{i+1}+\dots+S_n\\
  &=(I-T_i)\circ S_i+T_i\circ S_i+S_{i+1}+\dots+S_n
  =S_i+S_{i+1}+\dots+S_n.
}
Applying the inverse endomorphism $S$ side by side to this equality and again using the commuting property of the endomorphisms, it follows that
\Eq{TR}{
   T_i\circ R_{i-1}+(I-T_i)\circ R_{i+1}=R_i \qquad(i\in\{1,\dots,n\}).
}
Completely similarly we can also get that
\Eq{*}{
   t_i r_{i-1}+(1-t_i) r_{i+1}=r_i\qquad(i\in\{1,\dots,n\}).
}

To prove the $(R_k,r_k)$-convexity of $f$ for a fixed $k\in\{1,\dots,n\}$, let $x,y\in D$ be fixed and define 
\Eq{*}{
  u_i:=R_i(x)+(I-R_i)(y) \qquad(i\in\{0,\dots,n+1\}).
}
Then $u_0=x$ and $u_{n+1}=y$ hold and, due to the identities in \eq{TR}, it easily follows that
\Eq{*}{
  u_i=T_i u_{i-1}+(I-T_i)u_{i+1} \qquad(i\in\{1,\dots,n\}).
}
By the $(T_i,t_i)$-convexity of $f$, we have 
\Eq{ui}{
  f(u_i)\leq t_i f(u_{i-1})+(1-t_i)f(u_{i+1})\qquad(i\in\{1,\dots,n\}).
}
Define now, for $i\in\{0,\dots,n+1\}$, the coefficients $c_i$ as follows:
\Eq{ci}{
  c_i
  :=\begin{cases}
    0 &\mbox{if } i\in\{0,n+1\},\\
    \dfrac{r_k(s_0+\dots+s_{i-1})}{t_is_{i-1}} &\mbox{if } i\in\{1,\dots,k\},\\[3mm]
    \dfrac{(1-r_k)(s_i+\dots +s_n)}{(1-t_i)s_i} \quad &\mbox{if } i\in\{k+1,\dots,n\}.
    \end{cases}
}
Observe that all these coefficients are positive. Using the second equality in \eq{Ss} for $i=k$, the coefficient $c_k$ possesses also the following form
\Eq{*}{
c_k=\dfrac{r_k(s_0+\dots+s_{k-1})}{t_ks_{k-1}}
=\dfrac{(s_k+\dots+s_n)(s_0+\dots+s_{k-1})}{(s_0+\dots+s_n)(1-t_k)s_k}
=\dfrac{(1-r_k)(s_k+\dots +s_n)}{(1-t_k)s_k},
}
that is, the second formula in \eq{ci} remains also valid for $i=k$.
In what follows, we show that these numbers satisfy the following system of linear equations:
\Eq{rci}{
  c_i&=(1-t_{i-1})c_{i-1}+t_{i+1}c_{i+1} \qquad
  \mbox{if}\quad i\in\{1,\dots,n\}\setminus\{k\}.
}
We prove this equality for $i\in\{1,\dots,k-1\}$ and for $i\in\{k+1,\dots,n\}$ separately. If $i\in\{1,\dots,k-1\}$, then $2\leq k$. Thus, 
\Eq{*}{
  c_1=\frac{r_ks_0}{t_1s_0}=\frac{r_k}{t_1} \qquad\mbox{and}\qquad 
  c_2=\frac{r_k(s_0+s_1)}{t_2s_1}=\frac{r_k(\frac{1-t_1}{t_1} s_1+s_1)}{t_2s_1}=\frac{r_k}{t_1t_2}.
}
Hence the equality $c_1=t_2c_2$ follows, which proves \eq{rci} for $i=1$. For $i\in\{2,\dots,k-1\}$, we have
\Eq{*}{
  (1-t_{i-1})c_{i-1}+t_{i+1}c_{i+1}
  &=(1-t_{i-1})\frac{r_k(s_0+\dots+s_{i-2})}{t_{i-1}s_{i-2}}
  +t_{i+1}\frac{r_k(s_0+\dots+s_{i})}{t_{i+1}s_{i}}\\
  &=t_i\frac{r_k(s_0+\dots+s_{i-2})}{t_is_{i-1}}
  +(1-t_{i})\dfrac{r_k(s_0+\dots+s_{i})}{t_{i}s_{i-1}}\\
  &=\frac{r_k(s_0+\dots+s_{i-2}+(1-t_{i})s_{i-1}+(1-t_{i})s_{i})}{t_is_{i-1}}=c_i. 
}
For the proof of \eq{rci} in the case $i\in\{k+1,\dots,n\}$, we consider first the subcase $i=n$. We now have $k\leq n-1$ and
\Eq{*}{
  c_n&=\frac{(1-r_k)s_n}{(1-t_n)s_n}=\frac{1-r_k}{1-t_n},\\
  c_{n-1}&=\frac{(1-r_k)(s_{n-1}+s_n)}{(1-t_{n-1})s_{n-1}}
  =\frac{(1-r_k)(s_{n-1}+\frac{t_n}{1-t_n}s_{n-1})}{(1-t_{n-1})s_{n-1}}
  =\frac{1-r_k}{(1-t_{n-1})(1-t_n)}.
}From here, we can see that $c_n=(1-t_{n-1})c_{n-1}$, which proves equality \eq{rci} in the case $i=n$. Finally, assume that $i\in\{k+1,\dots,n-1\}$. Then 
\Eq{*}{
  (1-t_{i-1})c_{i-1}+t_{i+1}c_{i+1}
  &=(1-t_{i-1})\frac{(1-r_k)(s_{i-1}+\dots+s_{n})}{(1-t_{i-1})s_{i-1}}
  +t_{i+1}\frac{(1-r_k)(s_{i+1}+\dots+s_{n})}{(1-t_{i+1})s_{i+1}}\\
  &=t_{i}\frac{(1-r_k)(s_{i-1}+\dots+s_{n})}{(1-t_{i})s_{i}}
  +(1-t_{i})\frac{(1-r_k)(s_{i+1}+\dots+s_{n})}{(1-t_{i})s_{i}}\\
  &=\frac{(1-r_k)(t_is_{i-1}+t_is_i+s_{i+1}+\dots+s_{n})}{(1-t_{i})s_{i}}=c_i. 
}
After these preparations, multiply the inequality \eq{ui} by $c_i$ side by side and sum up the resulting inequalities for $i\in\{1,\dots,n\}$. Then, in view of the equalities \eq{rci}, the terms containing $f(u_i)$ for $i\in\{1,\dots,k-1,k+1,\dots,n\}$ cancel out and we obtain
\Eq{*}{
  c_k f(u_k)
  \leq t_1c_1 f(u_0)+(1-t_{k-1})c_{k-1} f(u_k) + t_{k+1}c_{k+1} f(u_k) + (1-t_n)c_n f(u_{n+1}).
}
This is equivalent to the inequality 
\Eq{uk}{
  (c_k -(1-t_{k-1})c_{k-1}-t_{k+1}c_{k+1}) f(u_k) 
  \leq r_kf(x)+(1-r_k)f(y).
}
Observe that in each inequality of \eq{ui}, the sums of the coefficients on the left and right hand side are equal to each other. 
This remains valid after multiplying by $c_i$ and summing up the inequalities so obtained. In particular, this has to be true for the 
inequality \eq{uk}. As a result, it follows that $c_k -(1-t_{k-1})c_{k-1}-t_{k+1}c_{k+1}=1$. Hence, the inequality \eq{uk} proves that $f$ is $(R_k,r_k)$-convex.
\end{proof}

\Cor{Last1}{Let $(X,+,d)$ be a metric  Abelian group, let $n_0\in\N$ such that $\mu_d(n_0)>1$ and let $D$ be a closed bounded $n_0$-convex set. Let $f:D\to\R$, $n\in\N$ and let $(S_0,s_0),\dots,(S_n,s_n)\in\E^d(X)\times\,]0,1[\,$ be such that $S_0,\dots,S_n$ are pairwise commuting and that $S_{0}+S_1$, \dots, $S_{n-1}+S_n$ and $S:=S_0+\dots+S_n$  are bijections with inverses belonging to $\E^d(X)$ and denote $s:=s_0+\dots+s_n$. Assume that $f$ is $\big((S_{i-1}+S_i)^{-1}\circ S_i,(s_{i-1}+s_i)^{-1}s_i\big)$-convex for all $i\in\{1,\dots,n\}$. Then, for all $k\in\{1,\dots,n\}$, we have $\big(S^{-1}\circ(S_k+\cdots+S_n),s^{-1}(s_k+\dots+s_n)\big)\in\CR_f^d$.}

\begin{proof}
For $i\in\{1,\dots,n\}$, define
\Eq{*}{
  T_i:=(S_{i-1}+S_i)^{-1}\circ S_i \qquad \mbox{and}\qquad
  t_i:=(s_{i-1}+s_i)^{-1}s_i.
}
Then, these endomorphisms and constants satisfy all the assumptions of the previous theorem, furthermore, the equalities in \eq{S-s} are satisfied. Therefore, the conclusion of this result applies.
\end{proof}

The next corollary is a generalization of former results of Dar\'oczy--P\'ales \cite{DarPal87} and Kuhn \cite{Kuh84} which are related to the vector space setting.

\Cor{Last2}{Let $(X,+,d)$ be a metric  Abelian group, let $n_0\in\N$ such that $\mu_d(n_0)>1$ and let $D$ be a closed bounded $n_0$-convex set. Let $f:D\to\R$ be $(T,t)$-convex for some $(T,t)\in\E_d(X)\times\,]0,1[\,$, and let $n\in\N$ be such that $\pi_n$ is a bijection with  a $d$-bounded inverse. Then, for all $k\in\{1,\dots,n\}$, $f$ is also $(\pi_n^{-1}\circ \pi_k,n^{-1}k)$-convex.}

\begin{proof} Assume that $f$ is $(T,t)$-convex. In order to use the previous corollary, for all $i\in\{0,\dots,n-1\}$, define
\Eq{*}{
  S_{2i}:=T,\qquad 
  S_{2i+1}:=I-T,  \qquad \mbox{and}\qquad
  s_{2i}:=t,\qquad 
  s_{2i+1}:=1-t.
}
Then, for every $i\in\{1,\dots,2n-1\}$, we have that $S_{i-1}+S_i=I$, which obviously has a $d$-bounded inverse. We also have that $S=S_0+\dots+S_{2n-1}=n\cdot I=\pi_n$, which has a bounded inverse by our assumptions. Furthermore, $\big((S_{i-1}+S_i)^{-1}\circ S_i,(s_{i-1}+s_i)^{-1}s_i\big)$ is equal to $(T,t)$ for even $i$ and to $(I-T,1-t)$ for odd $i$, which shows that \cor{Last1} is applicable. Therefore, by the conclusion of this corollary, $f$ is $\big(S^{-1}\circ(S_{2(n-k)}+\cdots+S_{2n-1}),s^{-1}(s_{2(n-k)}+\dots+s_{2n-1})\big)$-convex, i.e., it is $(\pi_n^{-1}\circ \pi_k,n^{-1}k)$-convex.
\end{proof}


\providecommand{\bysame}{\leavevmode\hbox to3em{\hrulefill}\thinspace}
\providecommand{\MR}{\relax\ifhmode\unskip\space\fi MR }
\providecommand{\MRhref}[2]{%
  \href{http://www.ams.org/mathscinet-getitem?mr=#1}{#2}
}
\providecommand{\href}[2]{#2}

\end{document}